\documentclass[12pt]{amsart}
\usepackage{graphicx}
\begin{document}
\newtheorem{theorem}{Theorem}[section]
\newtheorem{prop}[theorem]{Proposition}
\newtheorem{lemma}[theorem]{Lemma}
\newtheorem{claim}[theorem]{Claim}
\newtheorem{cor}[theorem]{Corollary}
\newtheorem{defin}[theorem]{Definition}
\newtheorem{example}[theorem]{Example}
\newtheorem{conj}[theorem]{Conjecture}
\newcommand{\aaa}{\mbox{$\alpha$}}
\newcommand{\map}{\mbox{$\rightarrow$}}
\newcommand{\ccc}{\mbox{$\mathcal C$}}
\newcommand{\kkk}{\mbox{$\kappa$}}
\newcommand{\Aaa}{\mbox{$\mathcal A$}}
\newcommand{\Fff}{\mbox{$\mathcal F$}}  
\newcommand{\bbb}{\mbox{$\beta$}}
\newcommand{\mlp}{\mbox{$\mu^{+}_{l}$}}
\newcommand{\ml}{\mbox{$\mu_{l}$}}
\newcommand{\mr}{\mbox{$\mu_{r}$}}
\newcommand{\mlpm}{\mbox{$\mu_{l}^{\pm}$}}
\newcommand{\mrpm}{\mbox{$\mu_{r}^{\pm}$}}
\newcommand{\mlm}{\mbox{$\mu_{l}^{-}$}}
\newcommand{\mrp}{\mbox{$\mu_{r}^{+}$}}
\newcommand{\mrm}{\mbox{$\mu_{r}^{-}$}}
\newcommand{\mm}{\mbox{$\mu^-$}}
\newcommand{\mz}{\mbox{$\mu^{\bot}$}}
\newcommand{\omp}{\mbox{$0_{-}^{+}$}}
\newcommand{\opm}{\mbox{$0_{+}^{-}$}}
\newcommand{\opp}{\mbox{$0_{+}^{+}$}}
\newcommand{\rz}{\mbox{$\rho^{\bot}$}}
\newcommand{\rp}{\mbox{$\rho^+$}}
\newcommand{\rmm}{\mbox{$\rho^-$}}
\newcommand{\mpm}{\mbox{$\mu^{\pm}$}}
\newcommand{\mpp}{\mbox{$\mu^+$}}
\newcommand{\mt}{\mbox{$\mu^{t}$}}
\newcommand{\mb}{\mbox{$\mu_{b}$}}
\newcommand{\Pt}{\mbox{$\tilde{P}$}}
\newcommand{\rpq}{\mbox{$\rho^{(p,q)}$}}
\newcommand{\mpq}{\mbox{$\mu^{(p,q)}$}}
\newcommand{\sss}{\mbox{$\sigma$}} 
\newcommand{\Sss}{\mbox{$\mathcal S$}}  
\newcommand{\aub}{\mbox{$A \cup_P B$}}  
\newcommand{\xuy}{\mbox{$X \cup_Q Y$}}  
\newcommand{\huh}{\mbox{$H_1 \cup_F H_2$}}  
\newcommand{\Ss}{\mbox{$\Sigma$}}
\newcommand{\ddd}{\mbox{$\delta$}} 
\newcommand{\rrr}{\mbox{$\rho$}} 
\newcommand{\Ggg}{\mbox{$\Gamma$}}
\newcommand{\ggg}{\mbox{$\gamma$}}
\newcommand{\ttt}{\mbox{$\tau$}} 
\newcommand{\bdd}{\mbox{$\partial$}}
\newcommand{\zzz}{\mbox{$\zeta$}}
\newcommand{\Ddd}{\mbox{$\Delta$}}
\newcommand{\qb} {\mbox{$Q_B$}}
\newcommand{\lll}{\mbox{$\lambda$}}
\newcommand{\Lll}{\mbox{$\Lambda$}}
\newcommand{\inter}{\mbox{${\rm int}$}}

\theoremstyle{remark}
\newtheorem{remark}[theorem]{Remark}

\numberwithin{equation}{section}

\title[]{Unknotting tunnels and Seifert surfaces}

\author{Martin Scharlemann}
\address{\hskip-\parindent
        Mathematics Department\\
        University of California\\
        Santa Barbara, CA 93106\\
        USA}
\email{mgscharl@math.ucsb.edu}

\author{Abigail Thompson}
\address{\hskip-\parindent
        Mathematics Department\\
        University of California\\
        Davis, CA 95616\\
        USA}
\email{thompson@math.ucdavis.edu}

\date{\today} 
\thanks{Research supported in part by grants from the National Science 
Foundation}

\begin{abstract} 
Let $K$ be a knot with an unknotting tunnel $\gamma$ and suppose that 
$K$ is not a $2$-bridge knot.  There is an invariant $\rho = p/q \in 
\mathbb{Q}/2 \mathbb{Z}$, $p$ odd, defined for the pair $(K, \gamma)$.  

The invariant $\rho$ has interesting geometric properties:  It is often 
straightforward to calculate; e. g. for $K$ a torus knot and $\gamma$ an 
annulus-spanning arc, $\rho(K, \gamma) = 1$.  Although $\rho$ is 
defined abstractly, it is naturally revealed when $K \cup \gamma$ is 
put in thin position.  If $\rho \neq 1$ then there is a minimal genus 
Seifert surface $F$ for $K$ such that the tunnel $\gamma$ can be slid 
and isotoped to lie on $F$.  One consequence: if $\rho(K, \gamma) \neq 
1$ then $genus(K) > 1$.  This confirms a conjecture of Goda and 
Teragaito for pairs $(K, \gamma)$ with $\rho(K, \gamma) \neq 1$.
\end{abstract}
\maketitle

\section{Introductory comments}

In \cite{GST} the following conjecture of Morimoto's was established: 
if a knot $K \subset S^3$ has a single unknotting tunnel $\ggg$, then 
$\ggg$ can be moved to be level with respect to the natural height 
function on $K$ given by a minimal bridge presentation of $K$.  The 
repeated theme of the proof is that by ``thinning'' the $1$-complex $K 
\cup \ggg$ one can simplify its presentation until the tunnel is 
either a level arc or a level circuit.

The present paper was originally motivated by two questions.  One was 
a rather specialized conjecture of Goda and Teragaito: must a 
hyperbolic knot which has both genus and tunnel number one necessarily 
be a $2$-bridge knot?  A second question was this: Once the thinning 
process used in the proof of \cite{GST} stops because the tunnel becomes 
level, can thin position arguments still tell us more? 

With respect to the second question, it turns out that there is an 
obstruction to further useful motion of $\ggg$ that can be expressed 
as an element $\rho \in \frac{Q}{2Z}$.  Surprisingly, further 
investigation showed that, so long as $K$ is {\em not} $2$-bridge, the 
obstruction $\rho$ can be defined in a way completely independent of 
thin position and thereby can be viewed as an invariant of the pair 
$(K, \gamma)$.  Moreover, this apparently new invariant has useful 
properties: It is not hard to calculate.  If $\rho \neq 1$, then the 
tunnel can be isotoped onto a minimal genus Seifert surface.  This, in 
combination with some work \cite{EU} of Eudave-Munoz and Uchida, 
verifies a conjecture of Goda and Teragaito \cite{GT} in the case in 
which $\rho(K, \gamma) \neq 1$.  If $\rho = 1$ and the tunnel is a 
level edge, then the tunnel can be moved so that instead of connecting 
two maxima (say) of $K$, it connects two minima and vice versa.  A 
future paper will expand on this observation, addressing the 
technically difficult and rather specialized case in which $\rho = 1$ 
with the goal of verifying the Goda-Teragaito conjecture in this final 
case.

\bigskip

Here are a few technical notes on conventions and notation used in the 
arguments that follow:
\begin{enumerate}
\item For $X \subset M$ a polyhedron, $\eta(X)$ will 
denote a closed regular neighborhood, whereas (abusing 
notation slightly) $M - \eta(X)$ will mean the closed complement of $\eta(X)$ 
in $M$.  

\bigskip

\item Pairs of curves in 
surfaces will typically be regarded as having been isotoped to minimize 
the number of points in which they intersect.  Only occasionally is care required with 
this convention.  For example, if a surface $S$ containing curves 
$\aaa$ and $\bbb$ is cut open along a circle $c$ and the remnants 
$\aaa - c$ and $\bbb - c$ are isotoped in $S - c$ to minimize their 
intersection (not necessarily fixing $\aaa \cap c$ or $\bbb \cap c$) 
then when $S$ is reassembled, new intersections are introduced because 
of twisting around $c$.  In most contexts this 
won't matter, since it is the absence of intersections that typically
complicates an argument.  

\bigskip

\item When put in thin position as in \cite{GST}, a $1$-complex 
$\Ggg$ in $S^3$ will typically be regarded as having first been made 
``generic'' with respect to the given height function on $S^3$; 
that is, all vertices will be of valence $3$, with two edges incident 
from above (resp.  below) and one from below (resp.  above).  At any 
height there will be at most one critical point or vertex.  This 
convention leads to the following semantic problem: A process which puts 
$\Ggg$ in thin position typically terminates when an edge of $\Ggg$ is 
made level.  Then $\Ggg$ is no longer generic, but can be made generic 
by a small perturbation in which the height function on the edge 
becomes monotonic.  To describe this situation we will sometimes say 
that the edge is a ``perturbed'' level 
edge.

\end{enumerate} 

\section{Unknotted handlebodies in $S^3$ and their splitting spheres}

Consider a standard genus two handlebody $H$ in $S^3$ and suppose 
$\mpp, \mm, \mt$ are three non-parallel, non-separating meridian disks 
for $H$, fixed throughout our discussion.  Let $\Ss$ denote the 
$4$-punctured sphere $\bdd H - (\mpp \cup \mm)$, with boundary 
components $\mlpm, \mrpm$.  Let $\mz$ denote a fixed separating 
meridian disk of $H$ that is disjoint from $\mpm$ and intersects $\mt$ 
in a single arc.  There's a natural projection of $\Ss$ to the 
rectangle $I \times I$ so that $\mz$ projects to a horizontal 
bisector, $\mt$ to a vertical bisector, the two copies $\mlp$ and 
$\mrp$ of $\mpp$ in $\bdd \Ss$ project near the points $\bdd I \times 
\{ 1 \}$ and the two copies copies $\mlm$ and $\mrm$ of $\mm$ in $\bdd 
\Ss$ project near the points $\bdd I \times \{ -1 \}$. 

A {\em complete} pair of arcs in $\Ss$ will be a pair of arcs whose 
boundary has one point on each boundary component of $\Ss$.  A 
complete pair of arcs $\lll_{\infty}$ disjoint from $\mt$ is said to 
have {\em infinite slope} and a complete pair of arcs $\lll_{0}$ that 
is disjoint from $\mz$ is said to have slope $0$.  The union 
$\lll_{\infty} \cup \lll_{0}$ divides $\Ss$ into two copies of $I 
\times I$, which we'll call the front face and the back face of $\Ss$.  
See Figure 1.  There is a natural correspondence between proper 
isotopy classes of complete pairs of arcs in $\Ss$ and the extended 
rationals $p/q \in {\mathbb Q} \cup \infty$.  Here $|p|$ is the number 
of times one of the pair intersects $\mz$, $|q|$ is the number of 
times it intersects $\mt$, and the fraction is positive (resp.  
negative) if the pair (when isotoped to have minimal intersection with 
$\lll_{\infty} \cup \lll_{0}$) is incident to the lower left corner 
$\mlm$ on the front (resp.  back) face.  Note that a complete pair of 
arcs in $\Ss$ for which one end of each arc lies on $\mm$ and the 
other on $\mpp$ corresponds to a rational $p/q$ with $p$ odd.  See 
Figure 2.  

\begin{figure}
\centering
\includegraphics[width=.6\textwidth]{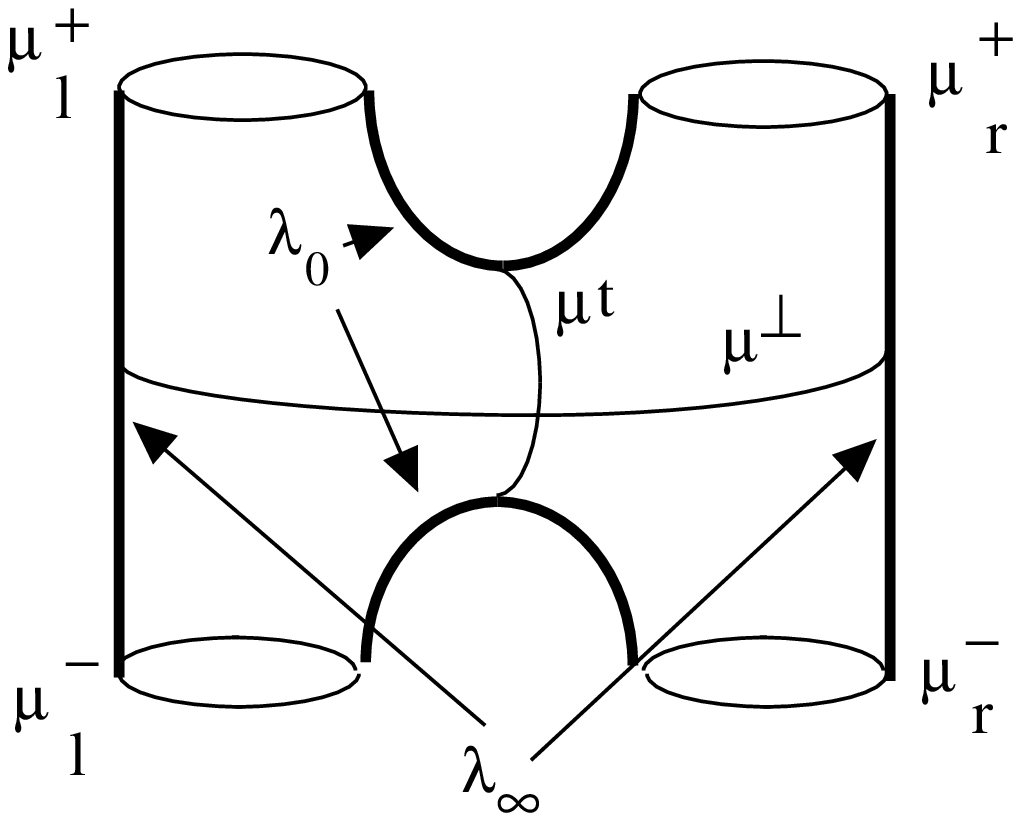}
\caption{}
\end{figure}  

\begin{defin} Given two complete pairs of arcs $\lll$ and $\lll'$, 
with slopes $p/q$ and $p'/q'$ respectively, define $\Delta(\lll, 
\lll') = |pq' - p'q|$.
\end{defin}

Note that if $\Delta(\lll, \lll') \leq 1$ then the two pairs can be 
isotoped to be disjoint; otherwise, $|\lll \cap \lll'| = 2\Delta(\lll, 
\lll') - 2$.  

\begin{figure}
\centering
\includegraphics[width=.6\textwidth]{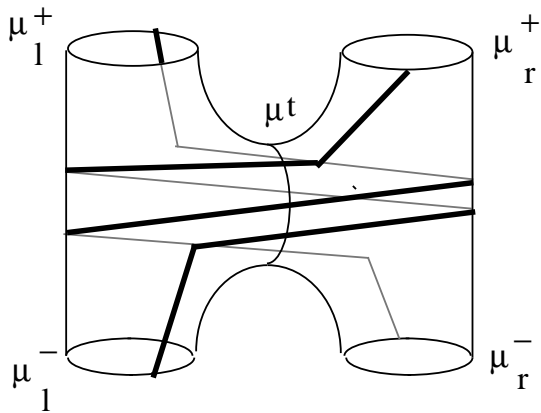}
\caption{Slope $1/3$}
\end{figure}  

\begin{defin} Suppose $H \subset S^3$ is an unknotted genus two 
handlebody and $S$ is a sphere that intersects $\bdd H$ in a single 
essential circle.  Then $S$ is a {\em splitting sphere} for $H$. 
\end{defin}

Alternatively, we could define a splitting sphere to be a reducing 
sphere for the Heegaard splitting $S^3 = H \cup_{\bdd H} (\overline{S^3 - H})$.

Suppose $S$ is a splitting sphere for $H \subset S^3$, so $D = S \cap 
H$ is an essential separating disk in $H$.  One possibility is that 
$D$ lies entirely inside the ball $H - (\mpp \cup \mm)$.  But, if not, 
then an outermost arc in $D$ of $D \cap (\mpp \cup \mm)$ cuts off a 
disk $D_0$ and $D_0 \cap \bdd H$ is an arc $\aaa$ with both ends at 
one of the boundary components of $\Ss$, say $\mlm$.  An arc in $\Ss$, 
such as $\aaa$, with both ends at $\mlm$, say, is called a {\em wave} 
based at $\mlm$.  A simple counting argument (as many ends of arcs 
of $\bdd D \cap \Ss$ lie on $\mlm$ as on $\mrm$) shows that there's 
also a wave $\aaa'$ based at $\mrm$ and that one of the components 
of $\Ss - \aaa$ is an annulus whose other end is one of $\mrp$ or 
$\mlp$.  A spanning arc for that annulus unambiguously gives us an arc 
with one end on $\mlm$ and the other end on one of $\mrp$ or $\mlp$.  
Similarly $\aaa'$ unambiguously gives us an arc from $\mrm$ to the 
other choice of $\mrp$ or $\mlp$.  Thus, given a splitting sphere, 
either its intersection circle with $\bdd H$ lies entirely in $\Ss$ or 
there is unambiguously defined a complete pair of arcs, each of which 
has one end on $\mm$ and one end on $\mpp$.  On the other hand, 
knowing that a specific essential pair $\lll$ is the result of this 
construction, we do not know whether the waves are based on $\mm$ or 
on $\mpp$.  In other words, to any choice of essential pairs of arcs, 
each of which has one end on $\mm$ and one end on $\mpp$, there 
correspond exactly two possible (pairs of) waves.

\begin{defin} Let $H$ be the standard genus two handlebody in $S^3$, 
$\mpp, \mm, \mt$ be three non-parallel, non-separating meridian disks 
for $H$ and $\mz$ a fourth, separating, meridian disk that is disjoint 
from $\mpm$ and which intersects $\mt$ in a single arc.  Finally let $S$ 
be a splitting sphere for $H$.

Define $\rho_{\bot}(\mpp, \mm, \mt, \mz, S) \in \mathbb{Q} \cup \infty$ 
to be the  slope associated to the waves of $S \cap \Ss$ as defined above.  

Two splitting spheres $S, S'$ are said to have {\em the same augmented 
slope} (with respect to $\mpp, \mm, \mt, \mz$), if $$\rho_{\bot}(\mpp, 
\mm, \mt, \mz, S) = \rho_{\bot}(\mpp, \mm, \mt, \mz, S')$$ and the 
associated waves are based at the same meridian $\mpp$ or $\mm$.
\end{defin}

A natural question is to what extent $\rho_{\bot}(\mpp, \mm, \mt,\mz, 
S)$, or indeed the augmented slope, depends on our choices.  Let's begin 
by considering different choices of splitting spheres.

\begin{defin} 
Let $S$ and $S'$ be two splitting spheres for $H 
\subset S^3$ and let $C$ and $C'$ be the corresponding separating 
curves $\bdd H \cap S$ and $\bdd H \cap S'$.  Define the {\em 
intersection number} $S \cdot S'$ to be the minimum number of points 
in $C \cap C'$.
\end{defin}

Since there is essentially only one way for a pair of circles $C$
and $C'$ in the sphere to intersect in either $2$ or $4$ points, the
relation between spheres with low intersection number is
easy to understand and describe.  Note, for example, that if $|C \cap
C'| = 4$ then each of the eight resulting arcs is adjacent to exactly
one bigon disk in each sphere.  See Figure 3.

\begin{figure}
\centering
\includegraphics[width=.4\textwidth]{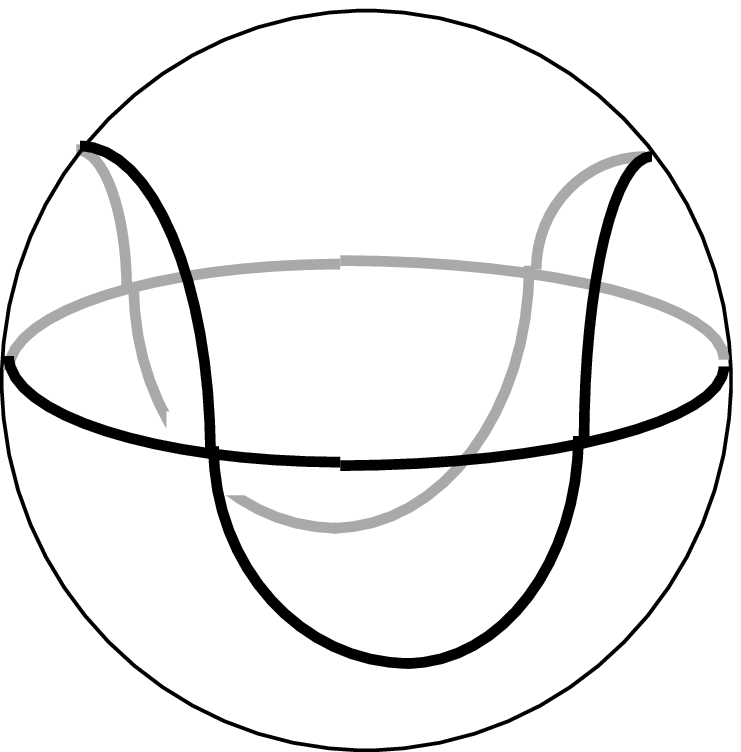}
\caption{}
\end{figure}

\begin{lemma}  \label{sphereintersect} If $S$ and $S'$ are not isotopic then
$S \cdot S' \geq 4$.   If $S \cdot S' = 4$ then $S
\cap S'$ is a single circle. Each of the two disjoint disks in $H$ obtained
by attaching a disk in $(S - S') \cap H$ to its adjacent bigon disk in
$(S' - S) \cap H$ is a non-separating essential disk in $H$ (and similarly for 
the handlebody $S^3 - H$). 
\end{lemma}

\begin{proof} If $S \cdot S' = 0$ then $S \cap S' = \emptyset$.  In a 
genus two handlebody, any pair of separating disks is parallel, so $S$ 
would be isotopic to $S'$.  If $S \cdot S' = 2$ then $C$ and $C'$ 
would be two separating circles in $\bdd H$ that intersect in two 
points, hence they would be isotopic, a contradiction.  Suppose $S 
\cdot S' = 4$.  Then each of the disks described in the lemma has 
boundary a bigon in $\bdd H$ that can be made disjoint from $C$ and 
$C'$.  If any boundary bigon were inessential it could be used to 
reduce $|C \cap C'|$.  So each bigon in $\bdd H$ (hence each disk in 
$H$) is essential.  
Furthermore, any essential circle in $\bdd H$ that is disjoint from 
$C$ and $C'$ is non-separating.  See Figure 4. \end{proof}
 
\begin{figure}
\centering
\includegraphics[width=.5\textwidth]{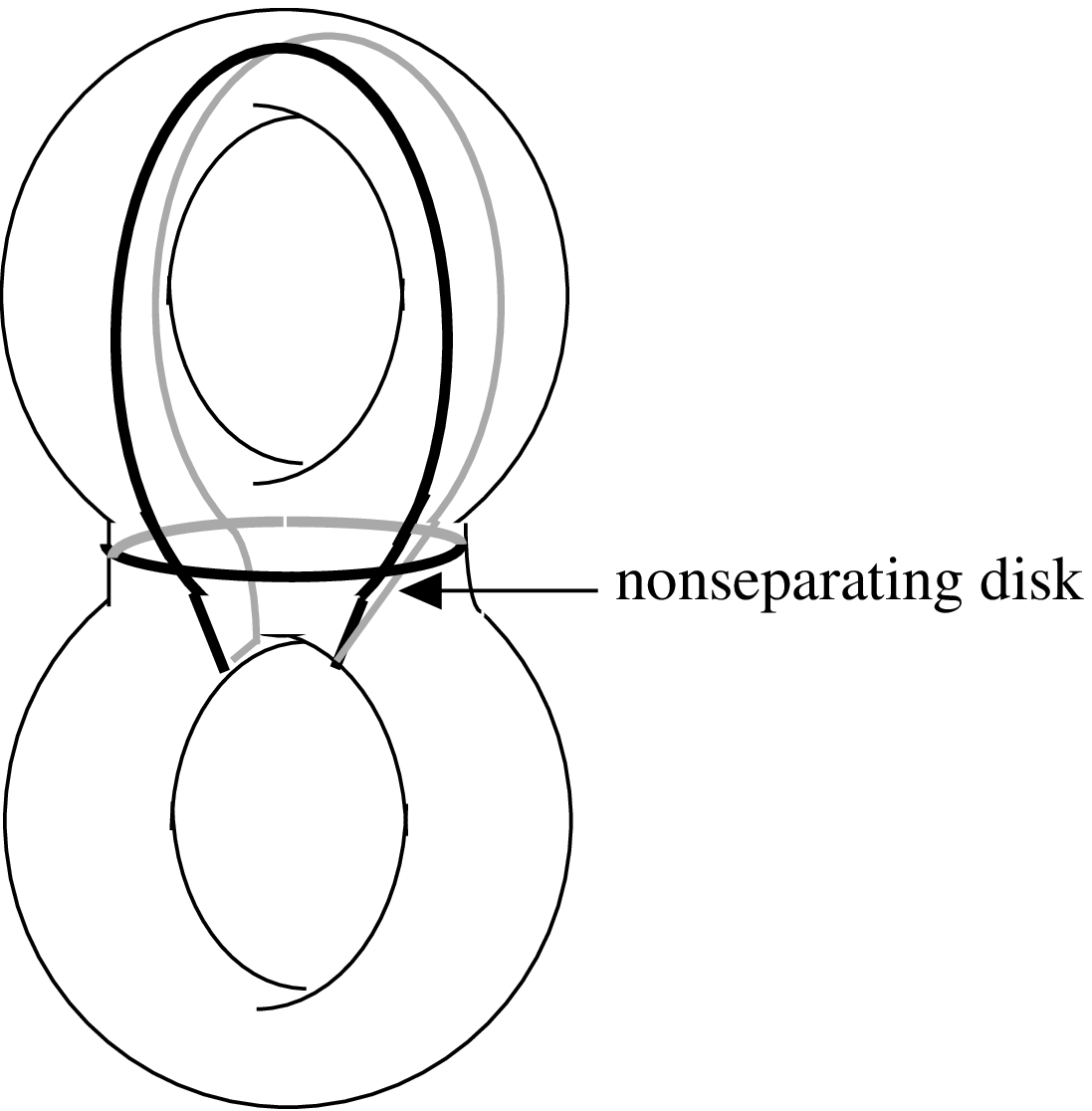}
\caption{}
\end{figure}  

\begin{prop}
If $S$ and $S'$ are two splitting spheres for $H \subset S^3$ then there is 
a sequence of splitting spheres $$S = S_0 \map S_1 \map \ldots \map S_m 
= S'$$ so that for $i = 1, \ldots, m$, $S_{i-1} \cdot S_i = 0$ or $4$.
\end{prop}

\begin{proof} There is an obvious (but obviously not unique) 
orientation-preserving homeomorphism $h: S^3 \map S^3$ with the 
property that $h(H) = H$ and $h(S) = S'$.  By the Alexander trick, $h$ 
is isotopic to the identity.

In  \cite{Go} Goeritz shows that any isotopy of $S^3$ that ends in a
homeomorphism carrying $H$ to $H$ is a product of particularly
simple such isotopies, whose effect on a fixed separating sphere $S_0$ 
is simple to describe.  In each case, either $S_0$ is preserved or the
intersection number of $S_0$ with its image is $4$.   

The upshot is this: the homeomorphism
$h$ is the composition of homeomorphisms $h = h_1 \circ h_{2} \circ
\ldots \circ h_m$ where each $h_i$ is the $H$-preserving homeomorphism
of $S^3$ obtained by one of the simple isotopies.  To obtain a sequence of
splitting spheres we  take $S_i = h_1 \circ h_2 \circ
\ldots \circ h_i(S_0)$.  Then notice that $S_i \cap S_{i-1}$
can be understood by viewing it as the image under the homeomorphism
$h_1 \circ h_2 \circ \ldots \circ h_{i-1}$ of $h_i(S_0) \cap
S_0$, so $S_i \cdot S_{i-1} = h_i(S_0) \cdot S_0 = 0$ or $4$.  \end{proof}   

{\bf Remark:}  This argument can be extended to Heegaard splittings of 
arbitrary genus, using work of Powell \cite{Po}.  See \cite{Sc}.

\begin{lemma} \label{knotcore} If two splitting spheres for $H$ have different
augmented slopes, then there is an essential disk in $S^3 -
H$ which intersects each of $\mpp$ and $\mm$ at most once. 
\end{lemma} 

\begin{proof}
The conclusion is obvious if any splitting sphere $S$ intersects $H$ 
in a disk disjoint from $\mpp \cup \mm$, for just use the disk $S - 
H$.  So we may as well assume that every splitting sphere defines an 
augmented slope.

Suppose $S$ and $S'$ are splitting spheres that give rise to two 
different augmented slopes.  Then there is a sequence of splitting 
spheres, beginning with $S$ and ending with $S'$, such that each has 
intersection number $4$ with the previous splitting sphere.  Since the first and 
last terms have different augmented slopes, somewhere there is a pair 
in sequence with different augmented slopes.  So we may as well assume 
that $S \cdot S' = 0$ or $4$.  

Let $\lll$ and $\lll'$ be the complete pair of 
arcs associated to the waves of $S \cap \bdd H$ and $S' \cap \bdd H$ 
respectively.  First notice that $\Delta(\lll, \lll') \leq 1$.  For if 
not, then $|\lll \cap \lll'| \geq 2$.  If we double each of the four 
arcs, the total number of intersection points is $8$, and converting a 
doubled arc into a wave can never remove intersection points, only add 
them.  

Suppose next that $\Delta(\lll, \lll') = 1$.  Then $\lll$ and $\lll'$ 
can be made disjoint, but not the waves that define them.  Indeed, if 
one pair of waves has its ends on $\mm$ and the other on $\mpp$ then 
each wave from $S$ intersects each wave from $S'$ in at least two 
points, a total of at least $2 \cdot 2 \cdot 2 = 8$ points.  See 
Figure 5.

\begin{figure}
\centering
\includegraphics[width=.4\textwidth]{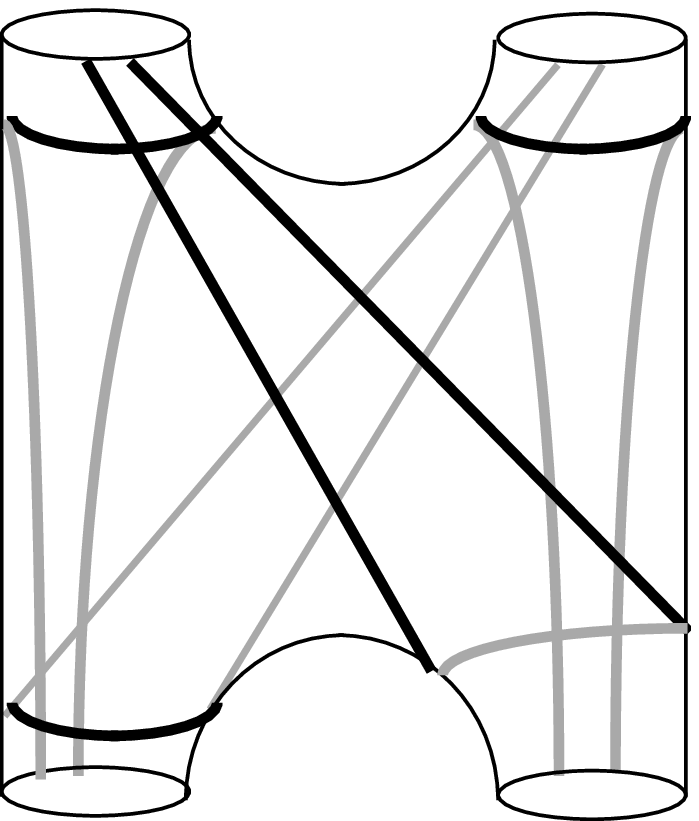}
\caption{}
\end{figure}

On the other hand, if $\Delta(\lll, \lll') = 1$ and both pairs of 
waves have their ends on $\mm$ (or both on $\mpp$) then each wave from 
$S$ intersects each pair of waves from $S'$ in at least $2$ points, a 
total intersection of {\em just the waves} of $4$ points.  Any other 
arc of $S$ with the same slope will intersect a wave of $S'$, and vice 
versa, so $S$ and $S'$ must both be disjoint from $\mpp$.  Following 
\ref{sphereintersect} we can say more: the bigons determined by the 
intersection points bound non-separating disks, two in $H$ and two in 
$S^3 - H$.  In our case each relevant bigon in $S^3 - H$ are made up 
of a unions of arcs, each with one end on $\mrm$ and one end on 
$\mlm$.  In particular, the bigon intersects $\mm$ twice or more, 
always with the same orientation.  But no such curve can bound a disk 
in $S^3 - H$, for the union of the solid torus $H - \mpp$ with the 
disk would define a punctured lens space in $S^3$.  See Figure 6.

\begin{figure}
\centering
\includegraphics[width=.6\textwidth]{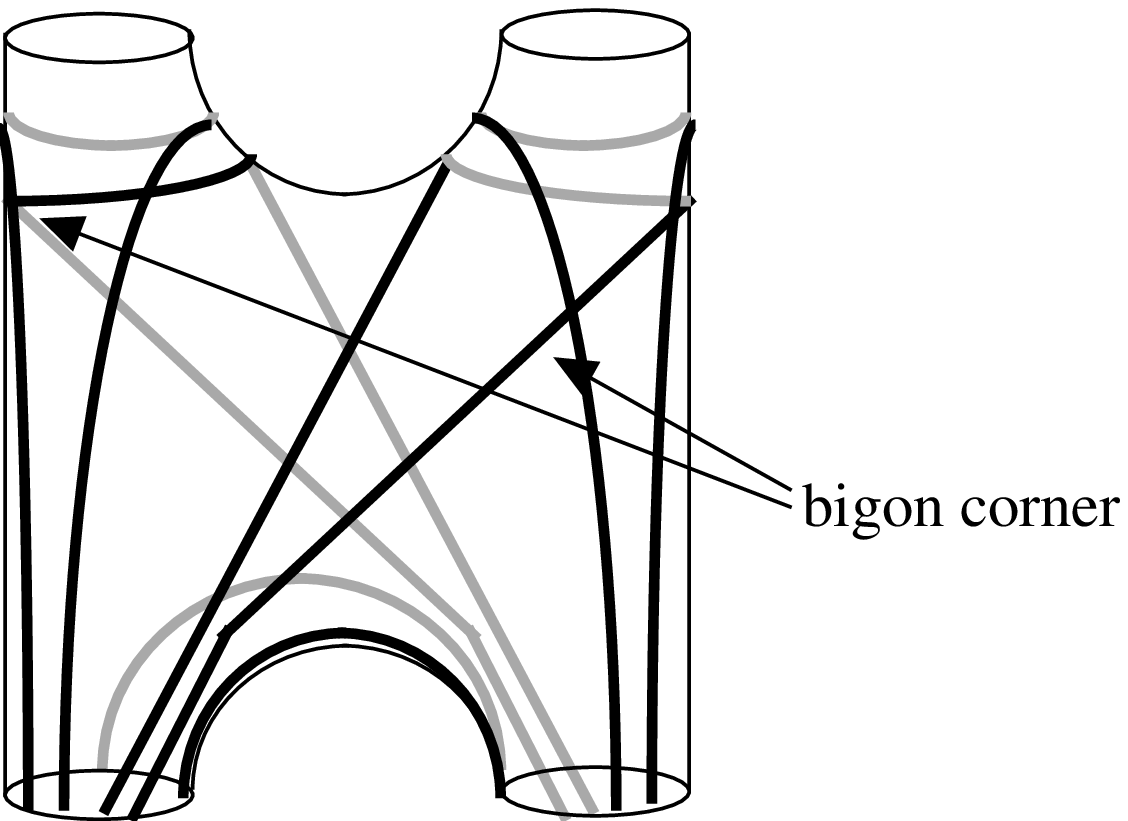}
\caption{}
\end{figure}

Finally, if $\Delta(\lll, \lll') = 0$ but $S$ has waves on $\mm$ and 
$S'$ on $\mpp$, then each wave of $S$ intersects each wave of $S'$ in 
at least (hence exactly) two points.  Because of the waves, any 
further component of $S \cap \Ss$ would intersect $S'$ and vice 
versa, so there can be no further components: $S$ and $S'$ each 
intersects $\Ss$ in exactly the two waves.  Then the bigons cut out 
by $C$ and $C'$ that bound disks in $S^3 - H$ each intersect each of 
$\mm$ and $\mpp$ in exactly one point.  \end{proof}

We have this corollary:

\begin{cor} \label{2bridge} If different splitting spheres define 
different augmented slopes for $(\mpp, \mm, \mt, \mz)$ then the knot 
core of the solid torus $H - \eta(\mt)$ is a $2$-bridge knot.
\end{cor}

\begin{proof} We can regard $H$ as the regular neighborhood of a 
$1$-vertex figure-$8$ graph $\Ggg$ in $S^3$, in which $\mpm$ are 
meridians of neighborhoods of to the two edges of the graph.  Let 
$k^{\pm} \subset \Ggg$ denote the subknots of $\Ggg$ corresponding to 
the meridian disks $\mpm$.  It suffices to show that $\Ggg$ is a 
standard unknotted figure-$8$ graph in $S^3$ since then the boundary 
of a regular neighborhood of the vertex of $\Ggg$ serves as a bridge 
sphere for a $2$-bridge presentation of the knot core of $H - 
\eta(\mt)$.  Following the unpublished \cite{HR} (see \cite{ST}) it 
suffices then to show that each of the knots $k^{\pm} \subset \Ggg$ is 
the unknot.

Suppose that $E \subset S^3 - H$ is an essential disk, given by Lemma 
\ref{knotcore}, that intersects each of $\bdd \mpp$ and $\bdd \mm$ at 
most once.  If $\bdd E$ is disjoint from exactly one of the meridians, 
say $\mm$, then $E$ is an unknotting disk for $k^+$, and $H \cup 
\eta(E)$ is an unknotted solid torus whose core is $k^-$.  Similarly, 
if $\bdd E$ is disjoint from both meridians, then $E$ 
divides $S^3 - H$ into two solid tori, each of whose meridians is an 
unknotting disk for one of $k^{\pm}$.  If $\bdd E$ intersects both 
meridians, then $H \cup \eta(E)$ is an unknotted solid torus in which 
both $k^+$ and $k^-$ can be viewed (individually) as core curves.  
\end{proof}

The next lemma shows that, given $\mt$, there is a natural choice of 
meridians $\mpm$.  First note that if $\rho_{\bot}(\mpp, \mm, \mt, \mz, 
S)$ is finite, then different choices of $\mz$ will change its value 
by a finite amount.  Indeed, any other possible $\mz$ will differ from 
the given one by some number of full Dehn twists around $\mt$, and 
such a Dehn twist changes $\rho_{\bot}$ by $\pm 2$.  

\begin{lemma} \label{madefinite} Suppose $\mt$ is a non-separating 
meridian disk for an unknotted handlebody $H$ and $S$ is a splitting 
sphere for $H$.  Then there is exactly one pair $\mpm$ of 
meridians disks for $H$ such that $\{ \mt, \mpp, \mm \}$ is a complete 
set of meridian disks for $H$ and $\rho_{\bot}(\mpp, \mm, \mt, \mz, S)$ is 
finite for any (hence every) choice of $\mz$ that is disjoint from $\mpm$ 
and intersects $\mt$ in a single arc.  
\end{lemma}

\begin{proof} Suppose, to begin, that there is an extension of $\mt$ 
to a set of meridians $\{ \mt, \mpp, \mm, \mz \}$ with respect to 
which $\rho_{\bot}$ is finite.  Because $\rho_{\bot}$ is finite, an 
outermost disk of $D = S \cap H$ cut off by the pair of meridians 
$\mpm \subset H$ intersects $\mt$.  Then an outermost subdisk $D_0$ of 
this subdisk, cut off by $\mt$, is disjoint from both meridians 
$\mpm$.  Furthermore, $\bdd D_0$ intersects $\mt$ in a single arc 
dividing $\mt$ into two subdisks.  The union of each of those subdisks 
with $D_0$ gives meridian disks for $H$ parallel to $\mpm$.  (See 
Figure 7).  

Now $\bdd D_0 - \mt$ is an essential arc $\aaa$ in the twice punctured 
torus $T_0 = \bdd H - \bdd \mt$, and the arc has both its ends on a 
single puncture.  Let $\aaa^{\pm}$ be closed curves in $T_0$ parallel 
to $\aaa \cup \{puncture \}$, lying on either side of $\aaa \cup 
\{puncture \}$.  If $\bbb$ is any other arc of $T_0 \cap S$ which has 
both its ends on a single puncture, then $\bbb$ is disjoint from 
$\aaa^{\pm}$; this is obvious if the ends of $\bbb$ lie on the other 
puncture, and follows from a counting argument on ends of arcs in $S 
\cap T_0$ if the ends of $\bbb$ lie on the same puncture as those of 
$\aaa$.  Any non-parallel separating pair of closed curves, e.  g.  
$\bdd \mpm$, in $T_0 - \aaa$ must be parallel to 
$\aaa^{\pm}$.  So we see that $\mpm$ are determined precisely by 
taking closed essential curves in $T_0$ that are parallel to $\aaa 
\cup \{puncture \}$.  (See Figure 8)

\begin{figure}
\centering
\includegraphics[width=.6\textwidth]{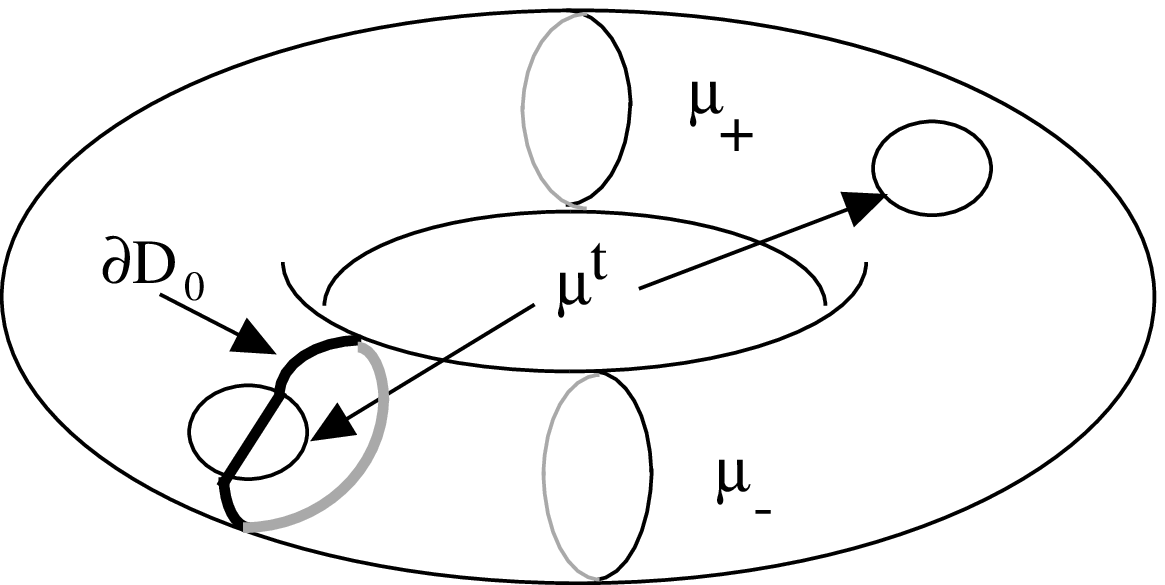}
\caption{}
\end{figure}

\begin{figure}
\centering
\includegraphics[width=.6\textwidth]{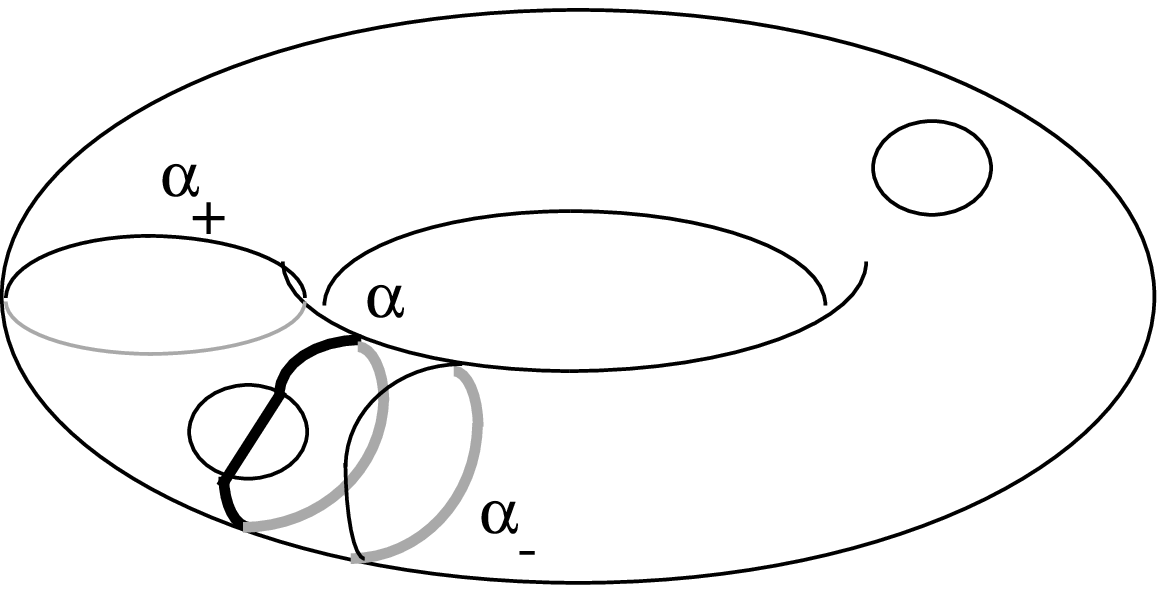}
\caption{}
\end{figure}

It's now easy to see that there is always some such pair.  Consider an 
outermost disk $D_0$ of $D$ cut off by $\mt$ in $H$.  Then the union 
of $D_0$ with each of the two subdisks of $\mt$ into which $D_0$ 
splits $\mt$ produces two natural meridian disks $\mpm$ for the solid 
torus $H - \mt$.  These, together 
with $\mt$ comprise a complete collection of meridian disks for $H$.  
Moreover, $\bdd D_0 \cap T_0$ is an essential arc that is disjoint 
from both meridians $\mpm$ of the solid torus $H - \mt$ bounded 
by $T_0$.  Since some arc of $S \cap T_0$ lying in the pairs of pants 
$T_0 - \mpm$ has both ends on $\bdd \mt$, no arc can have both ends on 
the same component of $\mpm$.  Hence $\rho_{\bot}$ is finite with 
respect to the meridian set $\{ \mpp, \mm, \mt \}$.
\end{proof}

\section{Knots with a single unknotting tunnel}

\begin{defin} A knot $K$ has {\em tunnel number one} if it is possible
to attach a single arc $\ggg$ to $K$ in $S^3 - K$ so that $S^3 - 
\eta(K \cup \ggg)$ is a solid handlebody.
\end{defin}

Put another way, $K$ has tunnel number one if there is a meridian disk
$\mt$ for a standard genus two handlebody $H \subset S^3$ so that the
solid torus $H - \mt$ has core knot isotopic to $K$.  

\begin{defin} Let $K$ be a knot and let $\ggg$ be an unknotting tunnel for 
$K$.  Let $S$ be a splitting sphere for $H = \eta(K \cup \ggg)$ and 
$\mt \subset H$ be a meridian disk for $\ggg$.  Let $\mpm$ be the 
meridians of $H - \mt$ given by Lemma \ref{madefinite} and $\mz$ be a 
meridian of $H - \mpm$ that intersects $\mt$ in a single arc.  Define 
$\rrr(K, \gamma, S) \in \mathbb{Q}/2\mathbb{Z}$ to be the value, mod 
$2$ of $\rrr_{\bot}(\mpp, \mm, \mt, \mz, S)$.
\end{defin}

Since different choices of $\mz$ change $\rrr_{\bot}$ by 
multiples of $2$, $\rrr(K, \ggg, S)$ is well-defined.  Moreover, by 
Corollary \ref{2bridge}, if $K$ is not 
$2$-bridge, $\rrr$ is independent of $S$ and so can be written 
$\rrr(K, \ggg)$.  Here we examine some features that $\rrr$ reveals 
about the knot and its tunnel.

Much is already known about their geometry.  The central theorem of 
\cite{GST} says that if the graph $K \cup \ggg$, viewed as a trivalent 
graph in $S^3$, is put in {\em thin position}, then $K$ is, on its 
own, in thin position and in bridge position, and $\ggg$ is a 
(perturbed) level edge.  Moreover, the tunnel $\ggg$ either has one 
end on each of two different maxima (or minima) or is a (perturbed) 
level loop, or ``eyeglass'', whose endpoint lies on a single maximum 
(or minimum) and which encircles all the other bridges of $K$.

The first claim is that the slope $\rrr(K, \ggg, S)$ is naturally 
revealed by some thin positioning of $K \cup \ggg$.  That is, there 
is a thin positioning of $K \cup \gamma$ so that the two isotopy 
classes of meridians of $K - \gamma$ with which level spheres 
intersect $K - \gamma$ are the classes $\mu^{\pm}$ identified in 
Lemma \ref{madefinite}.  

\bigskip

First we consider the case in which, upon thinning, $\ggg$ becomes a 
(perturbed) level eyeglass.  We can equivalently take, in this case, 
$\ggg$ to be a level edge with both ends incident to $K$ in the same 
point $p \in K$, i.  e.  $\ggg$ is a cycle.  Let $P$ be the level 
sphere in which $\ggg$ lies.  Following \cite{GST}, we can 
furthermore take $p$ to be the lowest maximum (or highest minimum) of 
the knot $K$; every bridge of $K$ other than the one containing 
$p$ has one end on each disk component of $P - \ggg$.  By a meridian 
of $K - p$ we will mean any meridian of $K$ that is disjoint from the 
(vertical) meridian of $K$ corresponding to the maximum $p$.

\begin{lemma} \label{wave1} Suppose $S$ is a splitting sphere for $H = 
\eta(K \cup\ggg)$ and let $D = S \cap H$.  Then any outermost disk 
$D_0$ of $D$ cut off by a meridian $\mu$ of $K - p$ intersects the 
meridian $\mt$ of the tunnel $\ggg$.  Moreover, a subdisk of $D_0$ cut 
off by $\mt$ is disjoint from a horizontal longitude of $\eta(\ggg)$.
\end{lemma}

\begin{proof} Cut $H$ open along the copies of $\mu$, denoted $\ml$ 
and $\mr$, corresponding to the points of $K \cap P$ that are nearest 
to $p$ in $K$.  Then $\mr$ and $\ml$ lie on different (disk) 
components of $P - \ggg$.  Cutting $H$ open along these meridians 
leaves one component that is a solid torus $W$ whose core is the cycle 
$\ggg$ and whose boundary contains disks corresponding to $\ml$ and 
$\mr$.  The meridian $\mt$ of $\ggg$ and the curves $P \cap \bdd W$ 
(parallel in $H$) naturally define respectively a meridian curve 
(which we continue to call $\mt$) and horizontal longitudes in the 
twice-punctured torus $T_0 = \bdd H \cap W$.  We want to understand 
the pattern of arcs $\Ggg = S \cap T_0$.  See Figure 9.

\begin{figure}
\centering
\includegraphics[width=.8\textwidth]{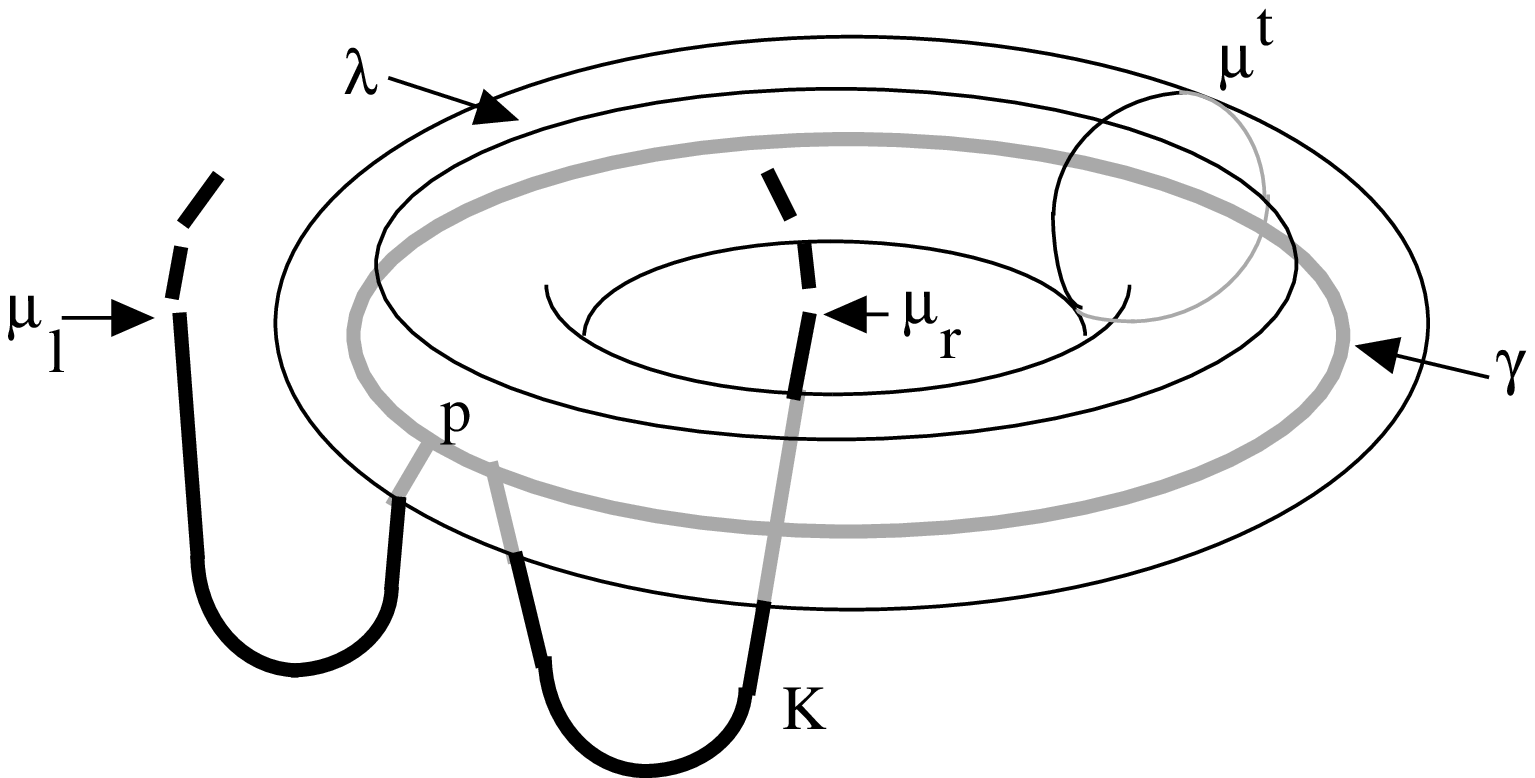}
\caption{}
\end{figure}

We begin by examining how $\Ggg$ intersects the twice punctured 
annulus $A$ obtained by cutting open $T_0$ along the horizontal 
longitude $\lll$ at the top of $\bdd H$. Note that $\mt$ intersects $A$ 
in a single spanning arc.  The boundary of $A$ can be thought of as 
two copies of $\lll$ which we denote $\bdd_0 A, \bdd_1 A$.

\bigskip

{\bf Claim}: Among the arcs in $\Ggg \cap A$ that intersect $\mt$, 
either there is one that has both ends on the same component of $\bdd 
A$ and separates the punctures $\mu_r$ and $\mu_l$, or there is one 
that has one end on $\bdd A$ and the 
other end on a puncture $\mu_r$ or $\mu_l$.

\bigskip

{\bf Proof of Claim:} Let $E = S - H$ be the exterior disk, and 
consider an outermost disk $E_0$ cut off from $E$ by an outermost arc 
$\aaa$ of $E \cap P$.  (It is easy to remove all closed components of 
$E \cap P$, since $K$ is thin.)  Since $\aaa$ obviously lies in a 
single component of $P - H$, it follows easily that $E_0$ must lie 
below $P$ (all arcs of $\Ggg - A$ can be assumed to be essential and 
so each spans the annulus $T_0 - A$) and so $E_0 \cap H$ lies in $A$.  
The arc $\aaa$ can't have one end on each of $\ml$ and $\mr$ since 
these lie in distinct components of $P - \ggg$.  If the ends of $\aaa$ 
both lie on the meridian $\ml$, say, then $\aaa$ is a longitudinal arc 
in the punctured annulus $A$.  That is, $\aaa \cup \ml$ is a core 
curve of $A$.  On the other hand, the outermost disk of $D$ cut off by 
$\mu$ must be a meridinal wave in $T_0$, so it too must also be based 
at $\ml$.  The complement of the two arcs, one meridinal and the other 
longitudinal, is then a disk in $T_0$ containing just the puncture 
$\mr$.  But then no wave could be based at $\mr$, and there would be 
more ends of $\Ggg$ on $\ml$ than on $\mr$, an impossibility.  We 
deduce that $\aaa$ has one or both ends on $\bdd A$.  Notice that if 
both ends of $\aaa$ lie on $\bdd A$ then they must lie on the same 
component ($\bdd_0 A$, say) of $\bdd A$ (since they are connected by 
an arc in $P - H$) and then the subdisk $A_0$ of $A$ cut off by $\aaa$ 
contains at least one puncture (or it would be inessential) but not 
both (else $\Ggg$ would intersect $\bdd_0 A$ more 
often than it intersects $\bdd_1 A$.

It remains to show that $\aaa$ intersects $\mt$.  We will show that if 
it doesn't, it can be used to make $K$ thinner.  Suppose that $\aaa$ 
is disjoint from $\mt$ and consider first the case in which $\aaa$ has 
one end on $\mr$ (say) and other end on $\bdd A$.  Then there is an 
arc $\kkk \subset (A - \mt)$ so that $\kkk \cup \aaa = \emptyset$ and 
one end of $\kkk$ lies on each of $\ml$ and $\mr$.  Since $\kkk$ is 
disjoint from $\mt$ it is isotopic in the ball $W - \mt$ to $K \cap (W 
- \mt)$.  Then $E_0$ can be used to pull a minimum of $K$ past $\kkk$ 
to $\aaa$, thinning $K$, a contradiction.  Similarly, if both ends of 
$\aaa$ lie on (necessarily the same component of) $\bdd A$, then $E_0$ 
is a lower cap, separating one component $K_-$ of $K - P$ below $P$ 
from all the others; $K_-$ is parallel to an arc $\kkk'$ in the 
punctured plane $P' = (P - H) \cup A$ that lies entirely in the twice 
punctured disk in $P'$ bounded by $\bdd E_0$.  There is then an arc 
$\kkk$ in $A - \mt$ whose interior is disjoint from $\kkk'$ and 
intersects $\aaa$ once and whose ends lie on the two punctures.  Then 
$K_-$ can be pushed up to $\kkk'$ and $K \cap W$ pushed down to 
$\kkk$, thinning $K$.  From this contradiction, we conclude that 
$\aaa$ intersects $\mt$, establishing the Claim.

\bigskip

Following the Claim, we have two cases to consider, corresponding to 
the two types of arcs given by the Claim.  Both arguments will 
use the $4$-punctured sphere $\Ss = T_0 - \mt$ bounded by $\ml, \mr$ 
and two copies $\mt_{\pm}$ of $\mt$.  We briefly recount some of its 
properties.  An outermost disk $D'$ of $D$ cut off by $\mt \cup \mu$ 
is a wave of $\Ggg$ in $\Ss$ and if it's based at $\ml$ an end count 
shows that there is a wave based at $\mr$ and vice versa.  Similarly 
if a wave is based at one of $\mt_{\pm}$ there is a wave based at the 
other.  Furthermore, any two waves 
in $\Ss$ must have the same slope.  The arc $\lll \cap \Ss$ is a path 
joining the two copies $\mt_{\pm}$; we can use it to establish slope 
$0 = 0/1$ in $\Ss$.  In these terms, a restatement of the Lemma is the 
claim that the wave determined by $D'$ is based at one of $\mt_{\pm}$ 
and is disjoint from $\lll$.

Suppose first that the arc $\aaa$ given by the Claim has both ends on 
a boundary component $\bdd_0 A$ of $A$ and cuts off from $A$ a disk 
$A_0$ containing a single puncture $\mr$ (say).  Once $\aaa \cap \mt$ is minimized 
by isotopy, an outermost arc of $\mt$ in $A_0$ cuts off a bigon $B$ 
containing $\mr$ and  bounded by subarcs of $\mt$ 
and $\aaa$. (See Figure 10.) An outermost arc of $\Ggg$ in $B$ (possibly the subarc of 
$\aaa$ in $\bdd B$) is a wave of $\Ggg$ in $\Ss$ based at $\mt_{\pm}$ 
that is disjoint from $\lll$, since $B$ is.  It follows that all the 
waves in $\Ss$, including that determined by $D'$, have the same 
property.

\begin{figure}
\centering
\includegraphics[width=.8\textwidth]{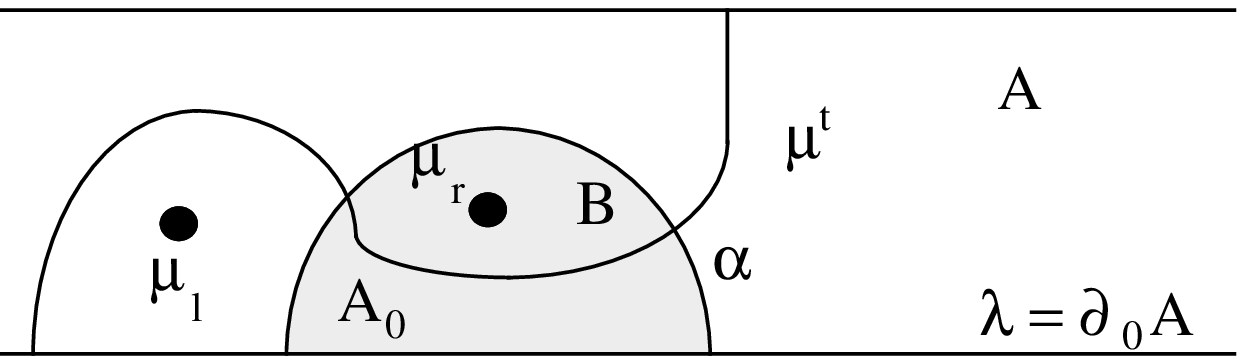}
\caption{}
\end{figure}

Suppose finally that only one end of $\aaa$ lies on $\lll$ and the 
other end lies on $\ml$ (say).  Since $\aaa$ intersects $\mt$, the end 
segments of $\aaa$ in $\Ss$ have these properties: One end, $\eta^+$, 
connects one of $\mt_{\pm}$ to $\ml$ and is disjoint from $\lll$.  The 
other end connects $\lll$ to one of $\mt_{\pm}$ essentially; let 
$\eta^-$ denote the segment of $\Ggg \cap \Ss$ that contains this end.  
Since one of $\eta^{\pm}$ intersects $\lll$ and one doesn't, they have 
different slopes in $\Ss$.  (See Figure 11.)  Observe that if waves 
are based on two boundary components of a $4$-punctured sphere, the 
only disjoint arc that can have a different slope than the waves is an 
arc that connects the bases of the waves.  Since only one end of 
$\eta^+$ can be the base of a wave, it follows that $\eta^-$ must 
connect the two bases of the waves.  Since one end of $\eta^-$ lies on 
one of $\mt_{\pm}$ this means that the waves must be based at 
$\mt_{\pm}$, and $\eta^-$ runs between $\mt_{\pm}$.  Finally, the 
slope of the wave must be that of $\eta^+$, so the wave is disjoint 
from $\lll$.
\end{proof}

\begin{figure}
\centering
\includegraphics[width=.4\textwidth]{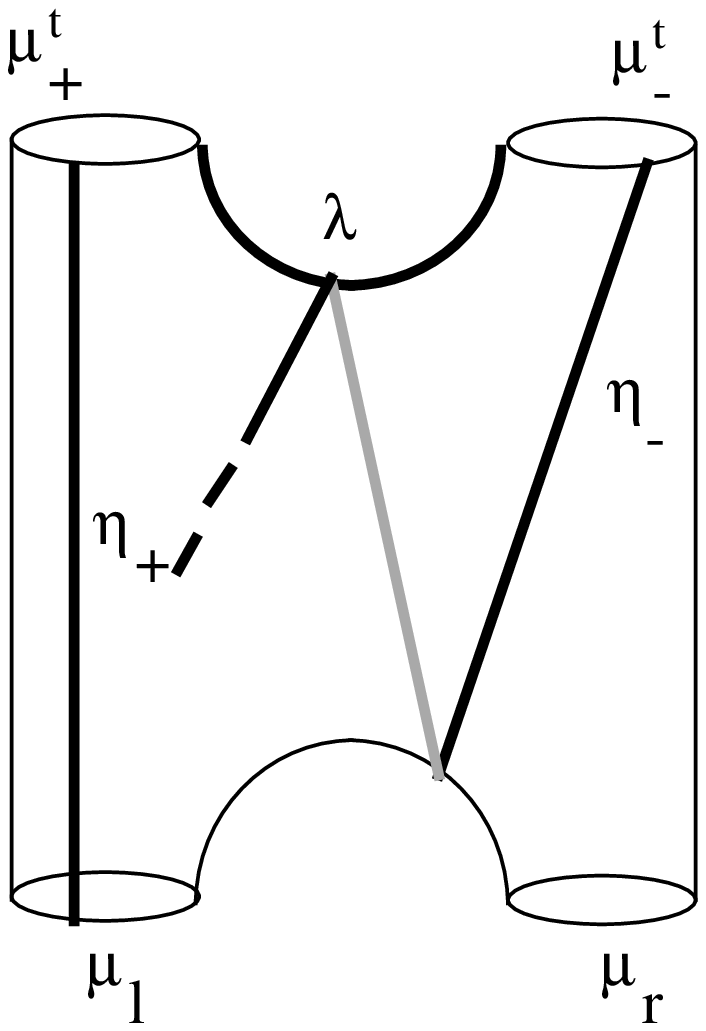}
\caption{}
\end{figure}  

\begin{cor} \label{wave1b} Suppose $S$ is a splitting sphere for $H = \eta(K 
\cup\ggg)$ and let $D = S \cap H$.  Suppose in a thin positioning of 
$K \cup \ggg$ the tunnel $\ggg$ is a level eyeglass at a maximum 
(minimum) $p$ of $K$.  Then the ends of $\ggg$ at $p$ may be slid 
slightly down (up) $K$ so that, with respect to the resulting 
meridians $\mpm$ of $K - \ggg$, $\rho$ is finite.
\end{cor}

\begin{proof} We know from Lemma \ref{madefinite} how to find 
meridians with respect to which $\rho$ is finite: Begin with an 
outermost disk $D_0$ of $D$ cut off by $\mt$ and choose meridians 
parallel to $\mt \cup \bdd D_0$ in the solid torus $H - \mt = 
\eta(K)$.  Lemma \ref{wave1} tells us precisely what those meridians 
are: one is $\mu$, the meridian of $K - p$.  The other is bounded by 
the union of an arc disjoint from the highest horizontal longitude 
$\lll$ of the circuit $\ggg_c$ and an arc that intersects $\lll$ once.  
Hence the boundary of the second meridian has a single maximum, so it 
can be viewed as simply the vertical meridian of $K$ at $p$, 
separating the ends of $\ggg$.
\end{proof}

\begin{theorem} \label{wave2}
 Suppose $S$ is a splitting sphere for $H = \eta(K \cup \ggg)$ and let 
 $D = S \cap H$.  Then $\ggg$ can be slid and isotoped to some thin 
 positioning of $K \cup \ggg$ so that $\rho$ is finite with respect to 
 the meridians of $K - \ggg$.
\end{theorem}

\begin{proof} Put $K \cup \ggg$ in thin position so that $\ggg$ can be 
levelled.  If $\ggg$ is an eyeglass, the result follows from Corollary 
\ref{wave1b}, so assume $\ggg$ is a level edge $e$.  Following 
\cite{GST} we can assume that $K$ is in minimal bridge position and 
the ends of $\ggg$ connect the two highest maxima of $K$.  If $P_u$ is 
a level sphere just below $e$, then the part of $S^3 - H$ lying above 
$P_u$ is just a collar $(P_u - H) \times I$, so we may as well assume 
that the exterior disk $E = S - H$ intersects this product region in 
bigons whose boundaries each consist of an arc in $P_u - H$ and an arc 
in $\bdd H - P_u$.

Let $\mpp$ and $\mm$ be the meridians of the two arcs $K - \ggg$; we 
know that $\mpp$, say, cuts off an outermost disk $D_0$ of $D = S \cap 
H$.  If $D_0$ intersects the meridian $\mt$ of the tunnel $\ggg$ we 
are done, so suppose that $D_0$ is disjoint from $\mt$.  Somewhere 
below $P_u$ and above the highest minimum of $K$, there is a generic 
level sphere $P$ which cuts off both an upper disk $E_u$ and a lower 
disk $E_l$ from $E$.  Since $E_u$ and $D_0$ can be made disjoint, it 
follows that $\bdd E_u$ crosses the meridian $\mt$ of the tunnel at 
most once.

We are now in a position to apply the argument of \cite[Theorem 
5.3]{GST}, though now in the context that the edge $e$ disjoint from 
$P$ is $\ggg$, not a subarc of $K$: There cannot be simultaneously an 
upper cap and a lower cap that have disjoint boundaries in $P$, or $K$ 
could be thinned.  If there is an upper cap and a disjoint lower disk 
or a lower cap and a disjoint upper disk, then, as in \cite{GST} we 
can find such a pair for which the interior of the disk is disjoint 
from $P$.  Then thin position implies that there cannot be 
simultaneously an upper cap and a lower disk and, if there is a lower 
cap which is disjoint from an upper disk (whose interior is now 
disjoint from $P$) then we can ensure that the boundary of the upper 
disk runs across the tunnel, hence exactly once across the tunnel.  
Similarly, if there is an upper disk which is disjoint from a lower 
disk then we can ensure that the interior of the upper disk is 
disjoint from $P$ and is disjoint either from a lower cap or a lower 
disk whose interior is also disjoint from $P$.  Furthermore, the 
boundary of the upper disk must be incident to the tunnel, hence run 
exactly once across the tunnel.  But then the upper and lower disks 
describe how to push the tunnel down to, or even below, the level of 
at least one minimum.  Following \cite[Corollary 6.2, Cases 1 and 
2]{GST} this implies (almost, see next paragraph) that the tunnel can 
be pushed down to connect two minima, which we may make to be the 
lowest minima.  Furthermore, this operation lowers by two the number 
of extrema of $K$ (in our case, the number of minima) found in the 
interior of the component of $K - \ggg$ containing $\mpp$, the base of 
the waves.

The argument then continues; it finally fails when there is only one 
extremum in the component of $K - \ggg$ containing $\mpp$.  That is, 
in the argument above, it may finally happen that the two maxima to 
which the ends of $\ggg$ are attached have only one minimum lying 
between them, and it's the component on which $\mpp$ is found.  In 
this case, the upper disk pushes $\ggg$ down to $P$, and then $\ggg$ 
together with the segment $\bbb \subset (K - P)$ containing the 
minimum form an unknotted cycle.  If any other minimum could be pushed 
up above this cycle then, again following \cite[Corollary 6.2]{GST}, 
all the other minima could, and the cycle would bound a disk disjoint 
from $K$, a contradiction.  We conclude that the lower disk or cap is 
in fact a lower disk that pushes $\bbb$ up to $P$, thereby forming a 
level cycle.  It can be used to slide $\ggg$ to form a level eyeglass, 
at which point we appeal to Lemma \ref{wave1}.
\end{proof}

\section{Pushing the tunnel off of a Seifert surface}

The next section will show that if $\rho(K, \ggg) \neq 1$ then 
$\ggg$ can be pushed onto a minimal genus Seifert surface for $K$.  
Ironically, the first step is to show that $\ggg$ can be pushed 
completely off of such a surface, which we will show in this section.

Let $K$ be a knot in $S^3$ and $F$ be a Seifert surface for $K$.  We say 
that $K$ is {\em parallel} in $S^3$ to an imbedded curve $c$ in $F$ if there 
is an annulus $A$ imbedded in $S^3$ such that $A \cap F = \bdd A = 
\bdd F \cup c = K \cup c$.

\begin{lemma} \label{minSeif1} Let $K$ be a knot in $S^3$ and $F$ be 
a minimal genus Seifert surface for $K$.  Suppose $K$ is parallel in 
$S^3$ to a 
curve $c$ in $F$.  Then $K = bdd F$ is parallel to $c$ in $F$.
\end{lemma}

\begin{proof} Let $A$ be the annulus giving the parallelism between 
$K$ and $c$.  Let $\eta(A)$ be a neighborhood of $A$ containing a 
neighborhood of $K$.  Since $A$ is an annulus, we can think of 
$\eta(A)$ as being a ribbon-like neighborhood of $K$ itself.  In the 
complement of $\eta(A)$, the remnant of $F$ is a possibly disconnected 
surface $\bar{F}$, with three (preferred) longitudinal boundaries on 
the boundary of $\eta(A)$.  If $\bar{F}$ is disconnected 
(corresponding to the case in which $c$ is separating) then one of the 
components of $\bar{F}$ is a Seifert surface for $K$.  Since it can't 
be of lower genus than $F$, the other component must be an annulus, 
defining a parallelism between $K$ and $c$ in $F$, as required.

Suppose $c$ is 
non-separating.  Then zero-framed surgery on $K$ yields a manifold $M$ 
and ``caps off"' $\bar{F}$; call the capped-off surface 
$\bar{F'}$.  A capped-off version $F'$ of the Seifert 
surface $F$ also imbeds in $M$ and $F'$ and $\bar{F'}$ 
represent the same homology class in $M$.  Since $genus(\bar{F'})$ 
is less than $genus(F')$, it follows from work of Gabai 
\cite[Corollary 8.3]{Ga} that genus $K$ is less than genus $F$, a 
contradiction.
\end{proof}  

\begin{prop} \label{tunneloff} Suppose $K$ is a knot and $\ggg$ is an 
unknotting tunnel for $K$.  Then $\ggg$ may be slid and isotoped until 
it is disjoint from some minimal genus Seifert surface for $K$.
\end{prop}

\begin{proof}
First choose a minimal genus Seifert surface $F$ and slide and isotope 
$\gamma$, doing both so as to minimize the number of points of 
intersection between $\gamma$ and $F$.  The slides and isotopies may 
leave $\ggg$ as either an edge or an eyeglass.  (In the latter case, 
let $\ggg_a$ be the edge in $\ggg$ and $\ggg_c$ be the circuit.)  We 
aim to show that $\gamma\cap F = \emptyset.$

Suppose to the contrary that after the slides and isotopies 
$\gamma\cap F$ is non-empty.  Let $E$ be an essential disk in the 
handlebody $S^3-\eta(K\cup\gamma)$ chosen to minimize the number $|E 
\cap F|$ of 
components in $E \cap F$.  $|E \cap F|>0$ for otherwise the 
incompressible $F$ would lie in a solid torus, namely (a component of) 
$S^3 - \eta(K \cup \gamma \cup E)$, and so be a disk.  Furthermore, 
since $F$ is incompressible, we can assume that ${E\cap F}$ consists 
entirely of arcs.  

Let $e$ be an outermost arc of ${E\cap F}$ in $E$, 
cutting off a subdisk $E_0$ of $E$.  The arc $e$ is essential in $F$, 
for otherwise we could find a different essential disk intersecting 
$F$ in fewer components.   Let $f=\partial(E_0)-e$, an arc 
in $\bdd \eta(K \cup \gamma)$ with each end either on the longitude $\bdd F 
\subset \bdd \eta(K)$ or a meridian disk of $\ggg$ corresponding to a 
point of $\ggg \cap F$.

\begin{itemize}

\item If any meridian of $\gamma$ is incident to exactly one end of $f$, 
then we can use $E_0$ to describe a simple isotopy of $\gamma$ which 
reduces the number of intersections between $\gamma$ and $F$.

\item If {\em no} meridian of $\ggg$ is incident to and end of $f$, then 
both ends of $f$ lie on $\bdd F \subset \bdd \eta (K)$.  If the 
interior of $f$ runs over $\ggg$ we are done, for $f$ is disjoint 
from $F$.  If the interior of $f$ lies entirely in $\bdd \eta(K)$ then
$E_0$ would be a boundary compressing disk for $F$ (since $e$ is 
essential), contradicting the  minimality of $genus(F)$.

\end{itemize}

The only remaining possibility is that both ends of $f$ lie on 
the same meridian of $\ggg$.  In this case, $e$ forms a loop in $F$ 
and the ends of $f$ adjacent to $e$ both run along the same subarc 
$\ggg_0$ of $\ggg$.  
Since $f$ is disjoint from $F$, $\ggg_0$ either terminates in 
$\bdd \eta(K)$ or $\ggg$ is an eyeglass and $\ggg_0$ terminates in 
the interior vertex of $\ggg$.  

If $\ggg_0$ terminates in an end of $\ggg$ in $\eta(K)$ then, since 
the interor of $f$ is disjoint from $F$, $f$ must intersect $\bdd 
\eta(K)$ in either an inessential arc in the torus or in a 
longitudinal arc.  The former case is impossible, since if the disk 
bounded by the inessential arc did not contain the other end of $\ggg$ 
then it could be isotoped away and $E \cap F$ reduced.  If the disk 
did contain the other end of $\ggg$, then $\bdd E$ would cross one end 
of $\ggg$ more often than the other, an impossibility.  It follows 
that $f$ intersects the torus $\bdd \eta(K)$ in a longitudinal arc.  Then 
$\eta(\ggg_0 \cup E_0)$ is a thickened annulus $A$, defining a 
parallelism in $S^3$ between $K$ and the loop $e$ on $F$.  By Lemma 
\ref{minSeif1} that means the loop $e$ is parallel to $\bdd F$.  
Substituting $A$ for the annulus between $e$ and $\bdd F$ in $F$ would 
create a Seifert surface for $K$ with fewer intersections with $\ggg$, 
a contradiction.

So $\ggg$ is an eyeglass and $\ggg_0$ terminates in the interior 
vertex of $\ggg$.  If, nonetheless, the interior of $f$ intersects 
$\bdd \eta(K)$ this means that $f$ traverses the edge $\ggg_a \subset 
\ggg$ so $\ggg_a$ is disjoint from $F$.  In that 
case, we can just repeat the argument above, absorbing $\eta(\ggg_a)$ 
into $\eta(K)$.  So we can assume that $f$ lies entirely on $\bdd 
\eta(\ggg)$.  Now the component $Q$ of $\bdd \eta(\ggg) - F$ on which 
$f$ lies is either a punctured torus (if $F$ is disjoint from 
$\ggg_c$) or a pair of pants.  In the former case, consider the 
Seifert surface $F'$ obtained from $F$ by removing the meridian disk 
$\mu_f$ of $F 
\cap \eta(\ggg)$ on which the ends of $f$ lie and substituting $Q$.  $F'$ is 
of one higher genus than $F$, and intersects $\ggg$ in one fewer point.  
Surgery to $F'$ using $E_0$ reclaims the minimal genus without 
introducing another point.  Thus we get a contradiction to our choice 
of $F$.  

If $Q$ is a pair of pants a similar argument works: Since $F$ is 
incompressible, it follows that the loop $e$ bounds a disk in $F$.  
Since $f$ is essential, that disk contains exactly one of the other 
two meridian disks (call it $\mu_e$) of $\ggg$ in $F$ that correspond 
to boundary components of $Q$.  Remove the meridian disks $\mu_f$ and 
$\mu_e$ from $F$ and attach instead an annulus that runs parallel to 
the subarc of $\ggg$ (containing the interior vertex) that has ends at 
$\mu_e$ and $\mu_f$.  This creates a Seifert surface $F'$ of genus one 
greater than $F$, but having one fewer intersection point with $\ggg$.  
Now do surgery on $F'$ using $E_0$, deriving the same contradiction as 
above.
\end{proof}

\section{Pushing the tunnel onto a Seifert surface}

In this section we will show that if $\rho(K, \ggg) \neq 1$ then 
$\ggg$ can be pushed onto a minimal genus Seifert surface for $K$.  
The first step was taken in the previous section: In general, $\ggg$ 
can be pushed completely off of some such surface.  A difficulty is 
that, after the slides used to push the tunnel off the Seifert 
surface, we no longer know that the resulting meridians of $K - \ggg$ 
are the ones by which we defined $\rho(K, \ggg)$ above.  Our strategy 
will be to use the meridians $\mpm$ by which we defined $\rho$, but at 
this cost: On the two punctured torus $\bdd \eta(K) - \ggg$ we no 
longer can assume that $\bdd F$ is a standard longitude.  All we know 
is that it is a curve that is isotopic in the {\em unpunctured} torus 
to the standard longitude.  The point of the following lemma, is that 
this situation is not a serious obstacle to further analysis.

\begin{lemma} \label{arctunnel} Suppose $K$ is a knot with unknotting 
tunnel $\ggg$, $H = \eta(K \cup \ggg)$, and $K_{0} \subset \bdd H$ is 
a curve in the twice-punctured torus $\bdd H - \eta(\ggg)$ which, in 
the unpunctured torus $\bdd \eta(K)$ is isotopic to a standard 
longitude.  Suppose $\aaa$ is an arc in $\bdd H$ so that $\aaa \cap 
K_0 = \bdd \aaa$ and $\aaa$ traverses $\eta(\ggg)$ once.  Then $\aaa$ 
is an unknotting tunnel for $K_0$ and the pair $(K_0, \aaa)$ is 
equivalent (by slides and isotopies) to the pair $(K, \ggg)$.
\end{lemma}

\begin{proof} Since $\aaa$ traverses $\ggg$ once, we can shrink 
$\aaa$, dragging along its end points in $K_0$ until $\aaa$ is just a 
spanning arc of the annulus $\bdd H \cap \eta(\ggg)$.  (At this point, 
we can identify $\aaa$ with $\ggg$ but we cannot yet identify $K_0$ 
with $K$.)  $K_0$ is a possibly complicated curve lying on $T = \bdd 
\eta(K)$ and $K_0$ is incident to $\aaa = \ggg$ at the ends of $\ggg$.

Using a collar $T \times I$ of $T$ in $\eta(K)$ isotope $K_0 \subset 
T$ until it is a standard longitude lying on the boundary of the 
smaller tubular neighborhood $\eta_- = \eta(K) - (T \times I)$ of the 
core $K$.  Extend the isotopy to an ambient isotopy of $T$, i.  e.  a 
self-homeomorphism of $T \times I \subset \eta(K)$.  This extends the 
ends of $\aaa$ as a $2$-braid through $T \times I$; call the extended arc 
$\aaa_+$.  The construction shows that $\eta_- \cup \eta(\aaa_+)$ is 
isotopic in $H$ to $\eta(K_0 \cup \aaa)$.  Now make the braid trivial 
by absorbing it into $\eta_-$.  (This translates into slides of the 
ends of $\aaa_+$ on $\bdd \eta_-$).  Afterwards, $H - (\eta_- \cup 
\eta(\aaa_+))$ is just a collar of $\bdd H$.
\end{proof}

\begin{theorem} \label{tunnelpush1} Suppose $K$ is a knot with 
unknotting tunnel $\ggg$ and $S$ is a splitting sphere for the 
handlebody $H = \eta(K \cup \gamma)$ with $\rrr(K, \ggg, S) \neq 1$.  
Suppose further that $F$ is an incompressible Seifert surface for the 
knot $K$ and that $F$ is disjoint from $\ggg$.  Then $\ggg$ can be 
slid and isotoped until it lies on $F$.
\end{theorem}

\begin{proof} Let $K_0$ be the copy $F \cap \bdd H$ of $K$ in $\bdd 
H$.  Consider the hemispheres $E = S - H$ and $D = S \cap H$.  By 
definition of $\rrr$ there are meridians $\mpp$ and $\mm$ for $\eta(K) 
\subset H$ that realize the slope $\rrr$.  Then, in particular, there 
are subarcs of $\bdd E = \bdd D$ that are waves based at one of these 
meridians.  (Warning: we know little about how these meridians 
intersect $K_0$.)  If the exterior disk $E$ were disjoint from $F$, 
then $F$ would lie in a solid torus obtained by compressing $H$ to the 
outside along $E$, contradicting the assumption that $F$ is 
incompressible.  So $F \cap E \neq \emptyset$.  We can isotope $\bdd 
E$ and $\bdd F$ to have minimal intersection and then remove any 
closed components of $F \cap E$ since $F$ is incompressible.

Let $E_0$ be an 
outermost disk of $E$ cut off by $F$.  Then $\bdd E_0$ consists of two 
arcs, $\aaa$ lying on $\bdd H$ and $\bbb$ lying on $F$.  We may assume 
that $\bbb$ is essential in $F$, for otherwise the subdisk of $F$ it 
cuts off, together with $E_0$, would again give an essential disk in 
$S^3 - H$ that is disjoint from $F$.  An important observation is that 
the ends of $\aaa$ lie on $\bdd F$ and are {\em incident to the same 
side of $F$}.  That is, if $\bdd F$ is normally oriented, the 
orientation points into (say) $\aaa$ at both ends of $\aaa$.

It's also true that $\aaa$ must cross the meridian $\mt$ of $\ggg$ at 
least once.  For otherwise, $E_0$ would give a $\bdd$-compression of 
$F$ to $\bdd \eta(K)$, contradicting the fact that $F$ is an 
incompressible (hence $\bdd$-incompressible) Seifert surface for $K$.  
If $\aaa$ crosses $\mt$ exactly once, then, following Lemma  
\ref{arctunnel}, $\aaa$ is equivalent to $\ggg$, so $E_0$ provides a way of 
isotoping $\ggg$ to $F$, completing the proof.  Hence it suffices to 
show:

\bigskip

{\bf Claim}: Any subarc of $\bdd E \cap \bdd H$ whose interior is disjoint 
from $\bdd F$ and whose ends lie on the same side of $\bdd F$ crosses $\mt$ 
at most once.

\bigskip

{\bf Proof of claim:} Continue to denote this subarc by $\aaa$.  Let 
$\Ss$ denote the four-punctured sphere obtained by cutting open $\bdd H$ 
along the meridians $\mpp$ and $\mm$ of $K$.

\bigskip

{\bf Case 1:} $\bdd F$ intersects $\Ss$ in loops as well as arcs.

\bigskip

Say the loops are based at $\mlp$ and $\mrp$ (see Figure 12).

\begin{figure}
\centering
\includegraphics[width=.6\textwidth]{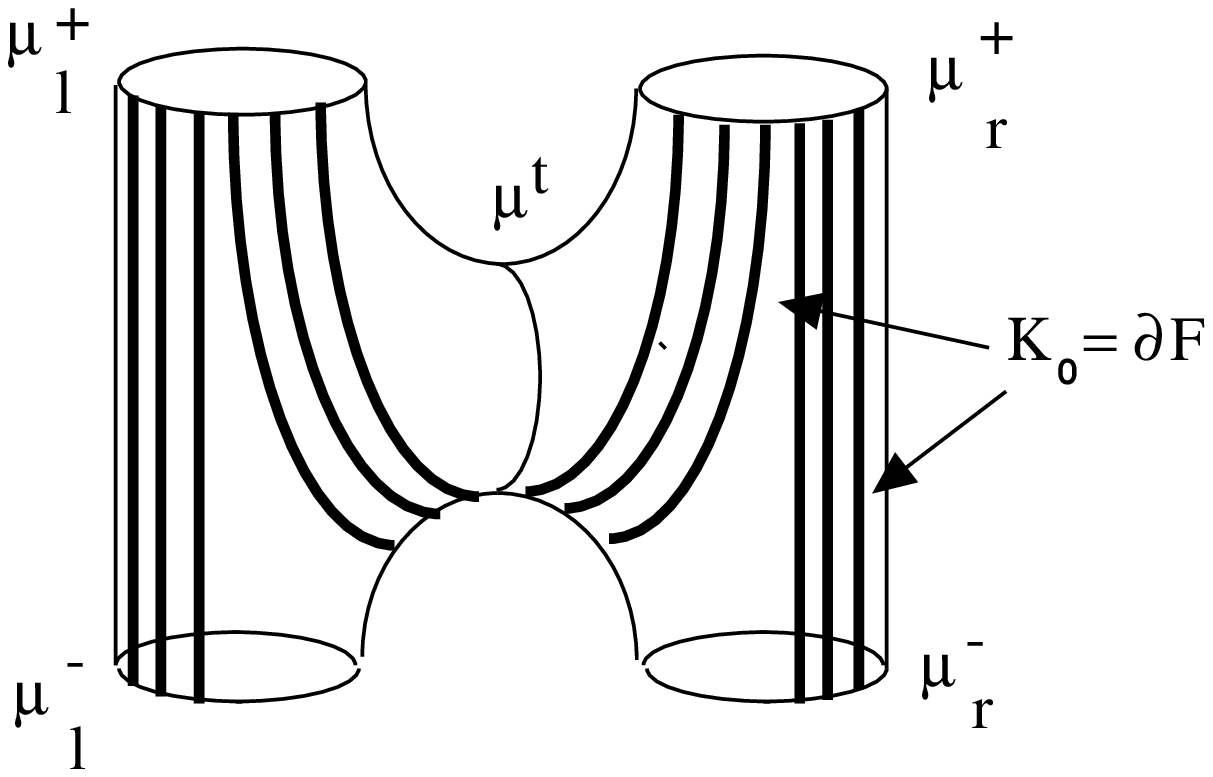}
\caption{}
\end{figure}

We first note that we may as well assume that $\aaa$ lies entirely in 
$\Ss$.  Indeed, since $\rrr$ is finite, there are waves of $\bdd E = 
\bdd D$ based at $\mt$ that lie in $\Ss$ and are on opposite sides of 
$\mt$ (see Figure 13).  It follows that $\aaa$ can't cross 
$\mm$.  For the same reason, we can assume that any component of $\aaa 
- \mpm$ with an end on $\mpp$ must lie in the component $\Ss^t$ of 
$\Ss - \bdd F$ that contains $\mt$, for any other component can be 
isotoped out through $\mpp$.  But the segments $\mrp \cap \Ss^t$ and 
$\mlp \cap \Ss^t$ are disjoint, for otherwise the two adjacent loops 
of $\bdd F \cap \Ss$ based at $\mlp$ and $\mrp$ would form a simple 
closed curve, which is impossible.

\begin{figure}
\centering
\includegraphics[width=.6\textwidth]{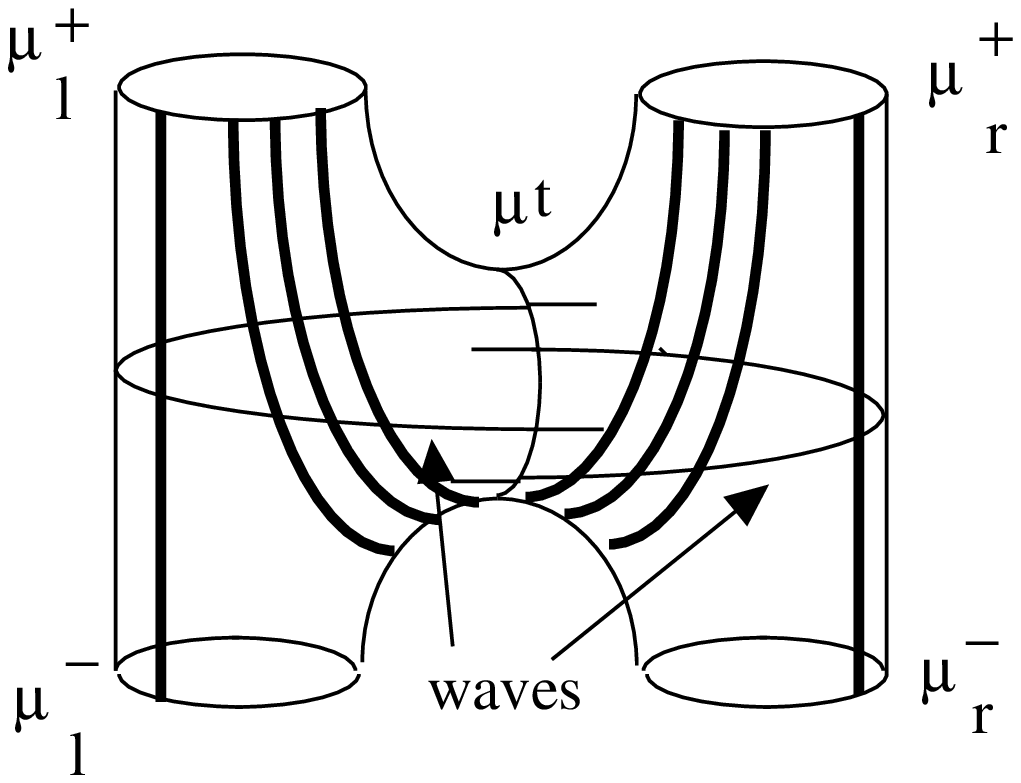}
\caption{}
\end{figure}

So now assume that $\aaa \subset \Ss$.  (We say $\aaa$ is {\em 
short}.)  Since $\bdd F$ intersects $\Ss$ in loops as well as arcs, 
then to intersect $\mt$ at all, $\aaa$ must lie in the annulus lying 
between the two outermost loops. (See Figure 14.) This annulus 
has $\mt$ as its core (since $\bdd F$ is disjoint from $\mt$) and any 
essential path in the annulus intersects $\mt$ at most once.

\begin{figure}
\centering
\includegraphics[width=.6\textwidth]{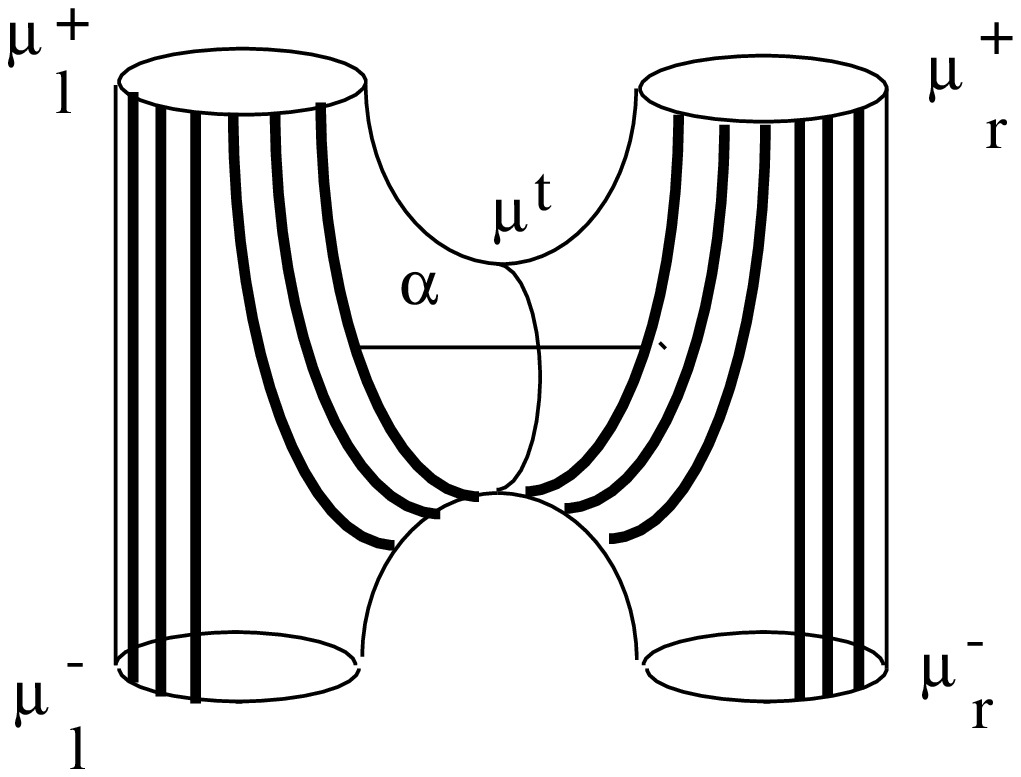}
\caption{}
\end{figure}

\bigskip

{\bf Case 2:} $\bdd F$ intersects $\Ss$ only in arcs.  

\bigskip

If $\bdd F$ intersects $\Ss$ only in arcs, then (considering the torus 
$H - \mt$) it must be in precisely two arcs, both of infinite slope 
(connecting $\mlp$ to $\mlm$ and $\mrp$ to $\mrm$.)  If $\aaa \subset 
\Ss$, the argument is the same as above, using the annulus $\Ss - \bdd 
F$ (see Figure 15).

\begin{figure}
\centering
\includegraphics[width=.6\textwidth]{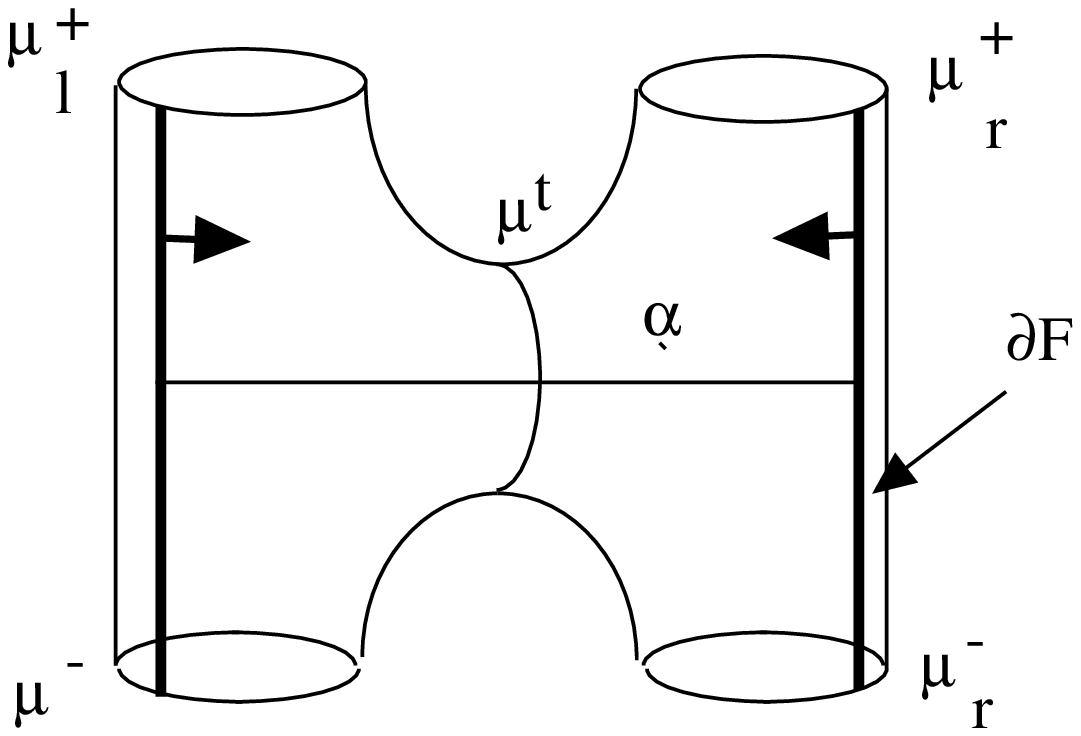}
\caption{}
\end{figure}

So finally suppose $\aaa$ is not contained in $\Ss$ and suppose with 
no loss that the waves of $\Ggg = S \cap \Ss$ are based at $\mlm$ and 
$\mrm$.  Any arc of $\Ggg = S \cap \Ss$ with an end on either $\mlp$ 
or $\mrp$ will then have a fixed slope, and since the slope $\rrr$ is 
not $+1 = -1 \in \frac{\mathbb{Q}}{2 \mathbb{Z}}$ it will intersect 
$\bdd F$ in its interior.  Moreover, if the normal orientation induced 
by that of $F$ points towards $\mlp$ on an arc with an end on $\mlp$ 
it will point away from $\mrp$ on an arc with an end on $\mrp$ (see 
Figure 16).  Hence we conclude that $\aaa$ cannot cross $\mpp$.  

It is as easy to rule out the possibility that $\aaa$ 
crosses $\mm$.  An arc of $\Ggg$ with one end on $\mlm$ may have slope 
$\rrr$ or have slope $\rho'$ with $\Delta(\rho, \rho') = 1$.  That is, 
fixing a meridian $\mz$ so that $\rho_{\bot} = p/q, q > |p|$, p odd, 
we could have that $\rho'_{\bot} = r/s$ with $|ps - rq| = 1$.  The arc couldn't 
have its other end on $\mpp$ so $r$ is even, hence $s$ is odd and 
$r/s \neq \pm 1$.  

Note that $|p/q - r/s| = 1/|qs|$ whereas both $|p/q|$ and $1 - |p/q| 
\geq 1/|q| \geq 1/|qs|$.  In words, $p/q$ is at least as close to 
$r/s$ as it is to $0$ or $\pm 1$.  Hence $|r/s| < 
1$ and either $r/s = 0$ (e.  g.  when $p/q = 1/k$) or $r/s$ has the 
same sign as $p/q$.  So any subarc of $\Ggg$ with an end on $\mlm$ or 
$\mrm$ is either disjoint from $\bdd F$ (if $\rho'_{\bot} = 0$) or it first 
intersects $\bdd F$ on the same side as an arc with slope $\rho$.  In 
particular, a subarc of $S \cap \bdd H$ that intersects $\mm$ and 
intersects $\bdd F$ precisely in its endpoints, necessarily ends on 
opposite sides of $\bdd F$, and so cannot be $\aaa$.
\end{proof}

\begin{figure}
\centering
\includegraphics[width=.6\textwidth]{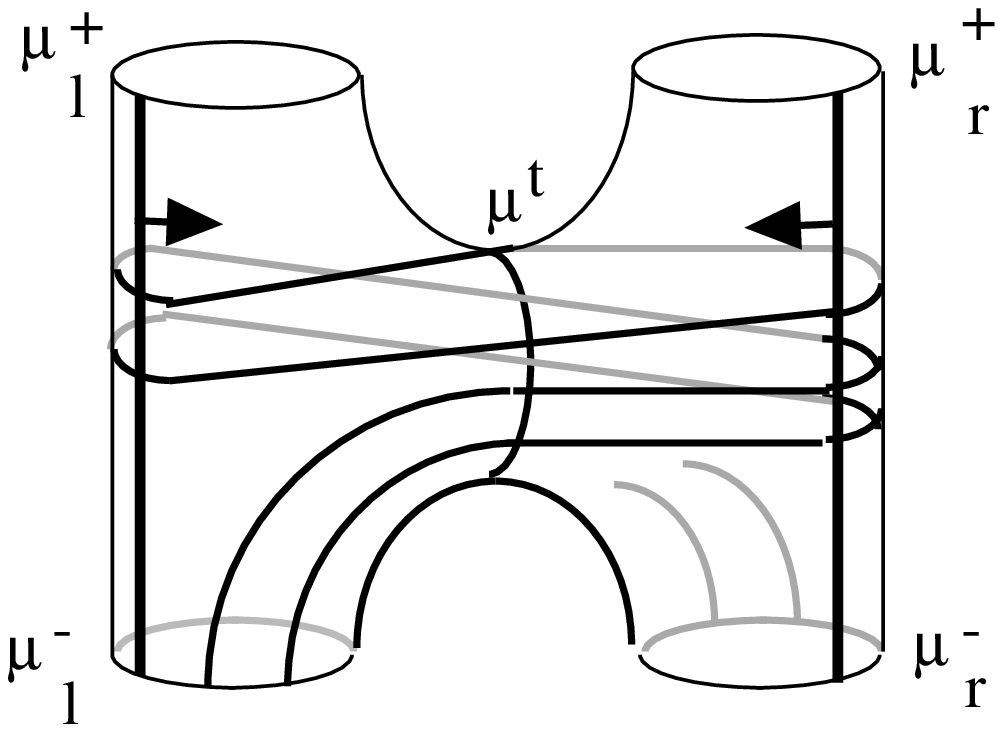}
\caption{}
\end{figure}

The proof of Case 1 in Theorem 
\ref{tunnelpush1}  did not require 
the assumption that $\rrr \neq 1$.  In particular, we have the 
following corollary:

\begin{cor} \label{tunnelpush2} Suppose $K$ is a knot with unknotting 
tunnel $\ggg$ and $S$ is a splitting sphere for the handlebody $H = 
\eta(K \cup \gamma)$.  Suppose further that $F$ is an incompressible 
Seifert surface for the knot $K$ and that $F$ is disjoint from $\ggg$.  
If either of the meridians $\mpm$ of $H$, chosen via \ref{madefinite} so that 
$\rrr$ is finite, intersects $K_0 = F \cap \bdd H$ in more than one 
point, then $\ggg$ is isotopic to an arc on $F$.
\end{cor}

More significantly: 

\begin{cor} \label{twobridge} Suppose $K$ is a knot with 
unknotting tunnel $\ggg$, $S$ is a splitting sphere for the handlebody 
$H = \eta(K \cup \gamma)$, $genus(K) = 1$ and $\rrr(K, \ggg, S) \neq 
1$.  Then $K$ is a $2$-bridge knot.
\end{cor}

\begin{proof} By Proposition \ref{tunneloff} we may assume that $\ggg$ 
is disjoint from some genus one Seifert surface $F$.  
Theorem \ref{tunnelpush1} shows that we can then isotope $\ggg$ onto 
$F$, necessarily as an essential arc.  Then $F - \eta(\ggg)$ is an 
incompressible annulus $A$ whose ends comprise a non-simple (because 
of $A$) tunnel number one link $L$.  (The core of $L$'s unknotting 
tunnel is the dual arc to 
$\ggg$ in the rectangle $\eta(\ggg) \cap F$.)  This implies, via \cite{EU}, that 
each component of $L$ is unknotted.  This then implies via \cite{HR} or 
\cite{ST} that the figure 8 graph obtained from $K \cup \ggg$ by 
crushing $\ggg$ to a point $v$ can be isotoped to lie in a plane.  This 
finally implies that $K$ is a $2$-bridge knot, with $\bdd \eta(v)$ the 
bridge sphere.  \end{proof}

This establishes the following conjecture of Goda-Teragaito (\cite{GT}) in the
case that $\rrr \neq 1$, and without the assumption that $K$ is 
hyperbolic.  

\begin{conj}[Goda-Teragaito]
A knot that is genus one, has tunnel number one, and is not a 
satellite knot is a $2$-bridge knot.
\end{conj}

The verification of the remaining case, when $\rho = 1$ and $K$ is 
hyperbolic, will be discussed elsewhere.  Note that Matsuda \cite{Ma} 
has verified the conjecture for all knots which are $1$-bridge on the 
uknotted torus, i. e. those with a $(1,1)$-decomposition.

\section{A sample calculation}

Let $T \subset S^3$ be an unknotted torus and $K \subset T$ be a torus 
knot in $T$.  Let $\ggg$ be a spanning arc for the annulus $T - K$.  
$\ggg$ is an unknotting tunnel for $K$ since $S^3 - \eta(K \cup \ggg)$ 
is a handlebody, namely the union of the interior and the exterior of 
$T$ along a disk in $T$.  In this section we will show that $\rho(K, 
\ggg) = 1$.  We will then use this calculation to construct examples 
of knots and tunnels with $\rho$ taking any value in $\mathbb{Q}/2 
\mathbb{Z}$.

To understand $\eta(K \cup \ggg)$ we will regard it as a bicollar of 
the punctured torus $T \cap \eta(K \cup \ggg) = \eta_{T}$ 
and consider its lift $\tilde{\eta} \times I$ to the universal cover $U = 
\mathbb{R}^2 \times I$ of $T \times I$.  Here's a back-handed way of 
doing that.  Since $\bar{L} = T - 
\eta_{T}$ is a disk, $\mathbb{R}^2 - \tilde{\eta}$ is a 
$\mathbb{Z}^2 = \mathbb{Z} \times \mathbb{Z}$ lattice of disks.  We 
can then regard $\tilde{\eta}$ as the complement of the slightly 
fattened lattice $L = \eta(\mathbb{Z}^2) \subset R^2$.

Using this picture, it is easy to describe how a lift of the meridian 
$\mt$ of the tunnel intersects $R^2$.  Begin by considering an arc 
$\mu$ connecting the lattice point $(0,0)$ with the lattice point $(m, 
n)$, where $m, n > 0$ are relatively prime.  Then the complement of 
all translates of $\mu$ in $\mathbb{R}^2 - L$ will be the complement 
of all translates of the line $my = nx$, namely an infinite collection 
of bands, each with slope $n/m$.  Each of these bands can also be 
described as a lift of $\eta_{T}(K)$ to $\mathbb{R}^2$, if $K$ is the 
$(n, m)$ torus knot.  Thus we see: if $K$ is the $(n, m)$ torus knot, 
$\mu \times I \subset \tilde{\eta} \times I$ is a lift of the meridian 
disk $\mt$ of $\ggg$ to $U$.

In a similar spirit, a vertical arc between adjacent points in 
$\mathbb{Z}^2$ is the lift of a meridian circle of $T$, and a 
horizontal arc between adjacent points is the lift of a longitude. 
Consider the following 
simple closed curve $\sss$ on $\bdd \eta(K \cup \ggg)$ (or rather the 
lift $\tilde{\sss}$ of $\sss$ to $\tilde{\eta} \times I$): $\tilde{\sss}$ intersects $R^2 \times \{ 0 \}$ 
in two adjacent vertical arcs connecting, say, the pair of points 
$(0,0)$ and $(0, 1)$ to the points $(1,0)$ and $(1,1)$.  
$\tilde{\sss}$ intersects $R^2 \times \{ 1 \}$ in the two
horizontal arcs connecting the pair of points $(0,0)$ and $(1,0)$ to 
the points $(0, 1)$ and $(1,1)$.  These two pairs of arcs, one 
vertical and one horizontal, are then 
connected to each other by product arcs in $\bdd L \times I$.  See 
Figure 17.  In particular, the curve $\tilde{\sss}$ projects to a 
unit square in $R^2$.

\begin{figure}
\centering
\includegraphics[width=.4\textwidth]{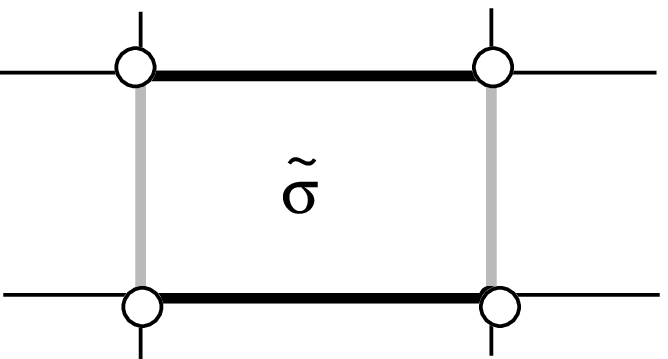}
\caption{}
\end{figure}

One can see that $\sss$ bounds an essential disk in both $\eta(K \cup 
\ggg)$ and in its complement.  Indeed, it bounds a disk $D$ in $\eta(K 
\cup \ggg)$ whose lift in $\tilde{\eta} \times I$ projects to the 
nullhomotopy of the square in the plane.  A disk $E$ that $\sss$ 
bounds in $S^3 - \eta(K \cup \ggg)$ can be described as the union of 
two meridian disks in each solid torus component of $S^3 - (T \times 
I)$ (a total of four disks) each attached along a single arc to a disk 
in the ball $\bar{L} \times I \subset T \times I$.  The disk $E \cap 
(\bar{L} \times I)$, when projected to $I$ has a single 
critical point, a saddle.

Now that we have found the tunnel meridian and a splitting sphere 
$S = D \cup E$, we 
need to find the preferred meridians of $K - \gamma$, that is, 
meridians with respect to which $\rho_{\bot}$ is finite.  The 
following argument is inspired by the proof of \cite[Lemma 2.2]{OZ}.  

Since $m, n$ are relatively prime, there are $p, q$ so that $0 < p < 
m$ and $0 < q < n$ and $mq - np = 1$.  Since $$det \left( 
\begin{array}{cc} m & n \\ p & q \end{array} \right) = 1,$$ it follows 
that $$\left( \begin{array}{cc} m & n \\ p & q \end{array} 
\right)^{-1}$$ is an integral matrix, so every point of $\mathbb{Z}^2$ 
is in $$\left( \begin{array}{cc} m & n \\ p & q \end{array} \right) 
\cdot \mathbb{Z}^2.$$ It follows that the parallelogram $P$ with 
corners $$\{ (0, 0), (m-p, n-q), (m, n), (p, q) \}$$ contains no 
element of $\mathbb{Z}^2$ in its interior.  See Figure 18.  In 
particular, appropriate lifts of the two triangles $$\Delta((0, 0), 
(m-p, n-q), (m, n))$$ and $$\Delta((0, 0), (p, q), (m, n))$$ tile each 
band in $\mathbb{R}^2$ that (as was described above) is a universal cover of 
$\eta_{T}(K) \subset T$.  It follows that the two arcs $\aaa^{\pm}$ in 
$\mathbb{R}^2$ with ends respectively at $(0, 0), (p, q)$ and $(0, 0), (m - p, 
m - q)$, when thickened in $U$, are lifts of meridians $\mpp$ and 
$\mm$ of $\eta(K)$.   It will now be straightforward 
to show that these are the appropriate meridians for calculating 
$\rho$.

\begin{figure}
\centering
\includegraphics[width=.6\textwidth]{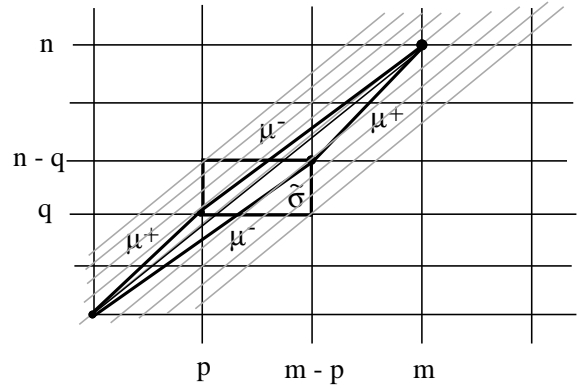}
\caption{The parallelogram P}
\end{figure}

We begin with the easy fact that if $a, b, c, d \geq 0$ are integers 
so that $a > b, c > d$ and $ac - bd = 1$, then $a = c = 1$ and $b = d 
= 0$.  It follows then from $det \left( \begin{array}{cc} m-p & n-q \\ 
p & q \end{array} \right) = 1$ that if $m - p > p$ then $n - q \geq q$ 
and if $m - p < p$ then $n - q \leq q$.  Rephrasing this as geometry: 
the minor diagonal of the parallelogram $P$, whose ends lie at $(p, 
q)$ and $(m-p, n-q)$, never has negative slope.  It follows that, for 
each corner $(p, q)$ and $(m-p, n-q)$ (say $(p, q)$) of $P$, some lift 
of the square $\sss$ has the property that it intersects $P$ in a 
triangle, with $(p, q)$ a vertex, and its entire opposite side lying 
in a single side of $P$.  (See Figure 19.)  Translating the plane 
geometry into the motivating context, the triangle represents a disk 
$D_0$, cut off from the interior disk $D$ bounded by $\sss$.  $D_0$ is 
cut off by the meridian of $K$ represented by the side of $P$ on which 
the triangle is based.  So the two sides of the triangle incident to 
$(p, q)$ are a wave of $\bdd D$ that crosses $\mt$, represented by the 
major diagonal of $P$, in two points.  Thus $\rho(K, \ggg)$ is not only seen 
to be finite, it is calculated to be $1 \in \mathbb{Q}/2 \mathbb{Z}$.

\begin{figure}
\centering
\includegraphics[width=.6\textwidth]{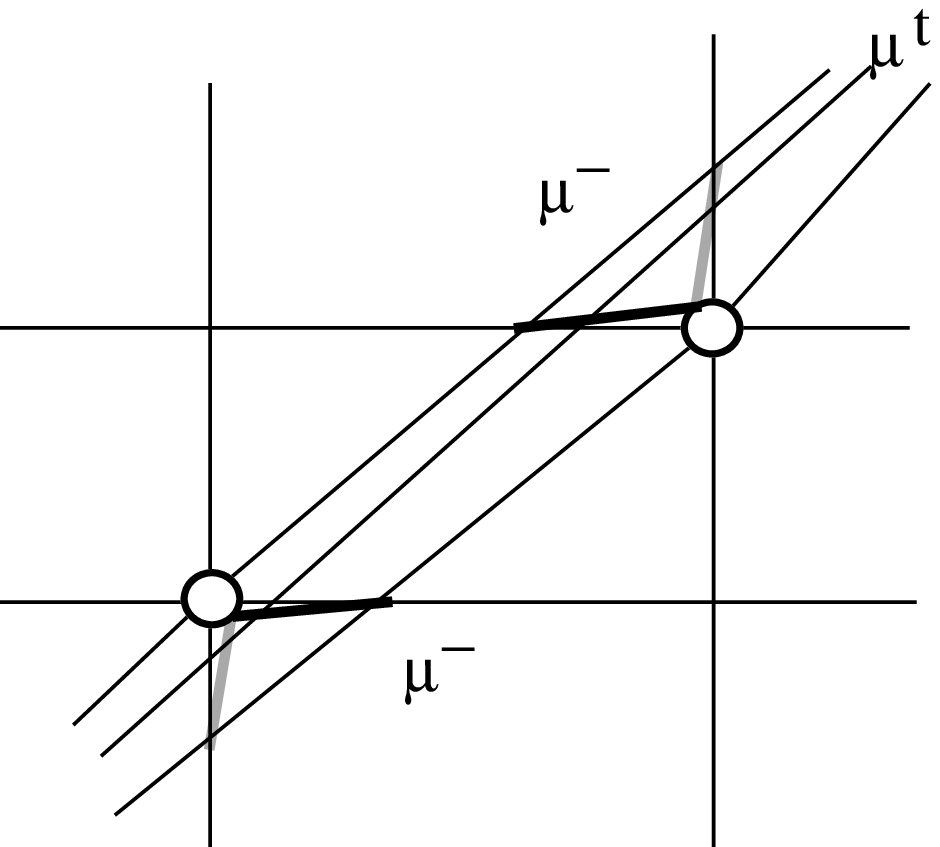}
\caption{}
\end{figure}

It is easy to extrapolate from this calculation to create examples in 
which $\rho$ can take on any value we like.  Namely, replace a 
the two strands of $K$ in a neighborhood of $\ggg$ by an appropriate 
rational tangle.  See Figure 20.  

\begin{figure}
\centering
\includegraphics[width=.6\textwidth]{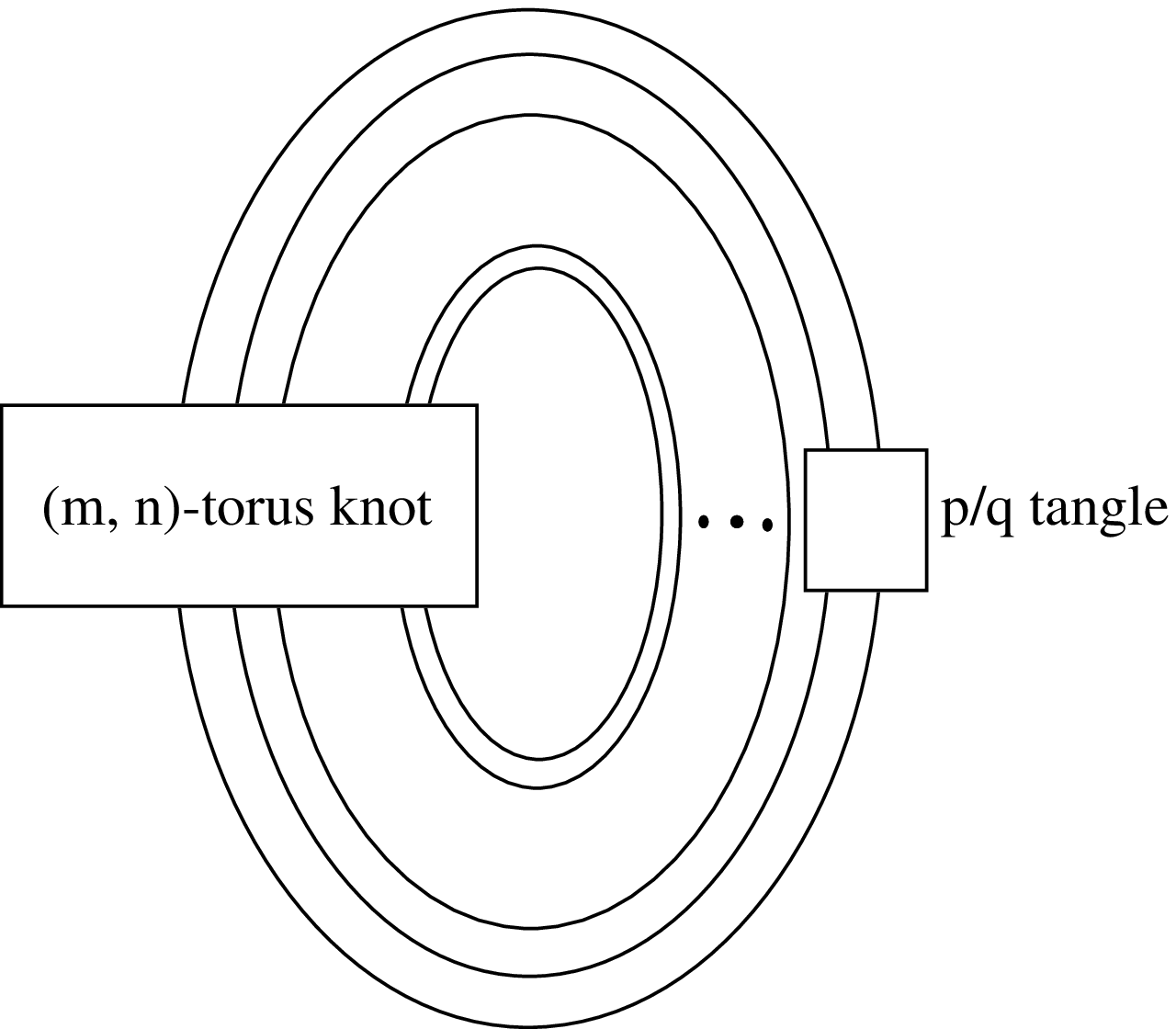}
\caption{}
\end{figure}

A torus knot is a simple example of a knot admitting a $(1, 
1)$-decomposition.  That is, it can be written as a $1$-bridge knot on 
an unknotted torus $T$in $S^3$ (see \cite{Do}, \cite{MS}).  Each knot 
$K$ admitting a $(1,1)$-decomposition has an unknotting tunnel $\ggg'$ 
(in fact two of them) best described as the eyeglass obtained by 
connecting the core $\ggg'_c$ of a solid torus $T$ bounds to the 
minimum of the knot $K \subset (T \times I)$ by an ascending arc 
$\ggg'_a$ in $T^2 \times I$.  (For 
$K$ a torus knot, $\ggg'$ is typically different from the tunnel 
$\ggg$ above.)  An ambitious reader should be able to discover an 
algorithm which determines, for this eyeglass tunnel $\ggg'$ and any knot 
$K$ admitting a $(1,1)$ decomposition, the value of $\rho(K, \ggg')$.


\end{document}